\documentclass[iicol,sn-basic]{sn-jnl}

\jyear{2022}%

\theoremstyle{thmstyleone}

\theoremstyle{thmstyletwo}

\theoremstyle{thmstylethree}

\newcommand{\design}{{\Omega_0}}
\newcommand{\density}{{\rho}}

\usepackage{caption}
\usepackage{subcaption}
\usepackage{booktabs}
\usepackage{dsfont}
\usepackage{algorithm,algpseudocode}
\usepackage{amsmath}

\usepackage{graphicx,booktabs,array}

\newcolumntype{M}[1]{>{\centering\arraybackslash}m{#1}}

\algdef{SE}[DOWHILE]{Do}{doWhile}{\algorithmicdo}[1]{\algorithmicwhile\ #1}

\begin{document}

\title[Article Title]{Topology Optimization with Frictional Self-Contact}

\author*[1]{\fnm{Zeshun} \sur{Zong}}\email{zeshunzong@math.ucla.edu}
\equalcont{These authors contributed equally to this work.}

\author[1]{\fnm{Xuan} \sur{Li}}\email{xuanli1@math.ucla.edu}
\equalcont{These authors contributed equally to this work.}

\author[2]{\fnm{Jianping} \sur{Ye}}\email{jpye00@umd.edu}
\author[1]{\fnm{Sian} \sur{Wen}}\email{sianwen22@g.ucla.edu}
\author[3]{\fnm{Yin} \sur{Yang}}\email{yin.yang@utah.edu}
\author[4]{\fnm{Danny M.} \sur{Kaufman}}\email{kaufman@adobe.com}
\author[1]{\fnm{Minchen} \sur{Li}}\email{minchen@math.ucla.edu}
\author*[1]{\fnm{Chenfanfu} \sur{Jiang}}\email{cffjiang@math.ucla.edu}

\affil[1]{\orgdiv{Department of Mathematics}, \orgname{UCLA}, \orgaddress{ \city{Los Angeles}, \state{CA}, \country{USA}}}

\affil[2]{\orgdiv{Department of Mathematics}, \orgname{University of Maryland}, \orgaddress{\city{College Park}, \state{MA}, \country{USA}}}

\affil[3]{\orgdiv{School of Computing}, \orgname{University of Utah}, \orgaddress{\city{Salt Lake City}, \state{UT}, \country{USA}}}

\affil[4]{\orgname{Adobe Research}, \orgaddress{\country{USA}}}

\abstract{Contact-aware topology optimization faces challenges in robustness, accuracy, and applicability to internal structural surfaces under self-contact. This work builds on the recently proposed barrier-based Incremental Potential Contact (IPC) model and presents a new self-contact-aware topology optimization framework. A combination of SIMP, adjoint sensitivity analysis, and the IPC frictional-contact model is investigated. Numerical examples for optimizing varying objective functions under contact are presented. The resulting algorithm proposed solves topology optimization for large deformation and complex frictionally contacting scenarios with accuracy and robustness.}


\maketitle

\begin{figure*}[h]
\centering
\includegraphics[width=0.9\textwidth]{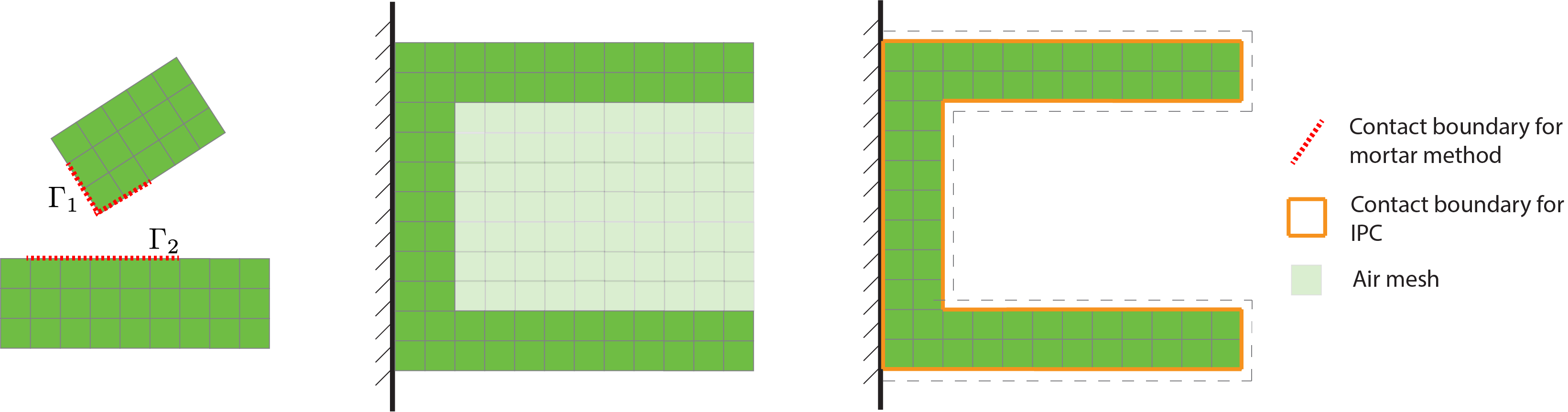}
\caption{Schematic graphical illustrations of three contact methods: mortar (left), air-mesh (mid), and IPC (right). Note that the contact boundary for the traditional mortar method must be manually specified and is typically a subset of the object boundary, while for IPC the contact boundary is the exact full boundary of the discrete mesh. The dashed line in the IPC graph is $\hat{d}-$away from the boundary of the solid structure, within which contact is treated.}\label{airmesh_v_ipc_scheme}
\end{figure*}

\section{Introduction}\label{sec1}
Topology optimization (TO) seeks to optimize material structural designs given user-specified inputs such as external loads and boundary conditions. It has been widely studied and developed for solving mechanical design problems across engineering fields \citep{sigmund_topoopt_approaches, SoftHAV}. 
Due to the often impractical assumption of small deformation, existing studies largely ignore contacting mechanisms. When large deformation is considered, however, the importance of dealing with contact in structural topology optimization becomes immediately apparent. In addition to preventing non-physical behaviors, i.e., interpenetrations, resolution of contact leads to different optimal structural designs that correctly take contact into account. Despite the need to accurately model contact to predict real-world behaviors, there has been little progress in optimizing topology with contact, especially when self-contact is required during deformation. The main challenges lie in (1) the lack of an accurate and robust model of contact that can be included in TO and (2) the appropriate resolution of complex contact behaviors, especially considering self-contact with friction \citep{sigmund_internal_contact}.

The literature on contact-aware topology optimization focuses primarily on two types of methods: mortar methods and fictitious domain methods. Fig. \ref{airmesh_v_ipc_scheme} shows a schematic illustration of both, as well as the IPC contact method we apply. 
Traditional mortar methods are generally limited to modeling the contact between a moving body and a fixed obstacle, often requiring a pre-specification (labeling) of contact surfaces. For each structural piece $i$, mortar methods pre-divide its boundary  $\partial \Omega_i$ into $\Gamma_i \subseteq \partial \Omega_i$ where a potential contact, and its complement $\Gamma_i^c \subseteq \partial \Omega_i$ where contact is not considered and so can not be resolved. Gap constraint functions are created, and a constrained optimization problem is solved. Satisfying the gap constraint ensures that nodes on the moving surface must not penetrate the element faces of the opposing obstacle surface \citep{hallquist1985sliding}, thus mimicking the contact between $\Gamma_i$'s. As contact interfaces require pre-specification, applications of mortar methods primarily focus on modeling either the interaction between elastic and fixed bodies as in \citep{kristiansen2020topology} or simple external contact such as \citep{Fernandez_mortar_2020, niu2019topology, mankame2004topology}. In addition, it is challenging for mortar methods to handle complex problems with nonlinear deformations and frictional effects. For example, \citet{Luo_mortar_hyperelastic} model contact with nonlinear springs for large deformations of hyperelastic bodies but cannot extend to frictional-contact cases, while \citet{Han_mortar} performs a node-to-node frictional analysis, but are limited to linear elasticity. 

Fictitious domain methods take a different path treating void regions between colliding bodies as a soft material with small stiffness. When two potentially colliding surfaces approach each other, the (filled) void region is compressed and so exerts large repulsion forces. The fictitious domain method was first introduced by \citet{Pagano_fictitious} to resolve self-contact. Many variations have followed.  \citet{Wriggers_third_medium} introduced the third medium contact method to handle external contact between bodies, which was then applied to self-contact \citep{sigmund_internal_contact, Muller_fictitious}. Despite the ability to resolve self-contact, these methods generally require manual hand-tuning of air-mesh parameters to avoid locking and parasitic transfer of non-physical forces. These methods also face significant challenges in modeling friction. Furthermore, the air mesh can potentially introduce large errors under significant distortions such as extreme shearing. Additional void regularization techniques are usually needed to alleviate such problems \citep{kruse2018isogeometric, Wriggers_third_medium}.

Recently, \citet{li2020incremental} propose a primal barrier-based Incremental Potential Contact (IPC) model for capturing the frictional contact of finite-strain elastic solids. Their method applies a smooth distance-based potential energy combined with a barrier-aware line search to avoid intersecting trajectories between surface primitives. Using a $C^2$ localized barrier, IPC automatically responds with contact forces between geometric pairs closer than a user-specified distance threshold and includes a corresponding variational friction model. As a result, IPC circumvents the difficulties covered above while providing guaranteed resolution of all contacting geometries. Note that while IPC was originally developed for elastodynamics, we have modified it here to solve for static force equilibrium in topology optimization. We assume hyperelasticity and present both mathematical and algorithm details for incorporating the IPC formulation into an existing topology optimization framework.

To summarize, we propose a new contact-aware topology optimization framework that can handle complex frictional contact scenarios, including external contacts and self-contacts under large deformation. A narrow-band process is adopted to define contact boundaries that evolve along the optimization procedure. An artificial timestep method is developed to properly model frictions within the static simulation scenario. Additional contributions include a strain limiting mechanism for tackling numerical difficulties introduced by low-density elements. We present numerical experiments to demonstrate the efficacy of the presented method. 

\section{Problem Statement}
Given a design domain $\design$, a standard density-based topology optimization problem seeks an optimal material distribution $\rho: \design\rightarrow [0,1]$ within $\design$ such that, under force equilibrium and a volume constraint, an applied objective function $F$ is minimized. Many approaches only consider external forces $f_{\text{ext}}$ and internal elastic forces $f_{\text{int}}$. Here we additionally include (normal-direction) contact forces and (tangential-direction) frictional forces.
Formally, the problem is
\begin{equation}
\begin{split}
    &\min_{\rho} F(\rho, u(\rho)) \;\; \text{s.t.} \\
    &\begin{cases} f_{\text{int}}(u)+ f_{\text{ext}}+ f_{\text{cont}}(u) + f_{\text{fric}}(u) = 0\\
    u\vert_{\Gamma} = u_0 \\
    V(\density) \leq \hat{V}, \end{cases}
\end{split} \label{eqn:problemstatement}
\end{equation}
where $u$ is the displacement field,  $\rho$ is the unknown scalar field describing the material allocation in $\design$, and $F(\rho, u(\rho))$ is the design objective function of interest. Here $f_{\text{cont}}$ and $f_{\text{fric}}$ are contact and friction forces applied on one or more pieces of structures in $\design$ due to their relative contact. A portion $\Gamma$ of the structure boundary $\partial \Omega_0$ will have prescribed Dirichlet boundary condition $u_0.$, while $V$ is the material's volume, $\int_{\design}\density dX$, that is constrained to be less than a user-specified upper bound criterion, $\hat{V}$. 

For actual manufacture, the material density $\rho$ should be close to either zero or one. Consequently, a binarization is conducted to the entries in $\rho$ for re-evaluating the objective function.

Problem (\ref{eqn:problemstatement}) applies to many widely studied topology optimization formulations. For instance, if $F$ is chosen to be the total elastic energy and  $f_{\text{cont}}=f_{\text{fric}}=0$, then we retrieve the classical structural compliance minimization problem \citep{sigmund88, sigmund99new, sigmund1997design}.

\section{Background}
\subsection{SIMP}
The Solid Isotropic Material with Penalization Method (SIMP) \citep{sigmund88, sigmund200199, sigmund99new} is widely applied in topology optimization. The local cell density $\rho_c \in [0,1]$ represents the material distribution on an Eulerian grid. SIMP assumes that Young's modulus of each grid cell $c$ is proportional to a polynomial of its cell density, i.e., $E_c = E_0 \rho_c^p$, where $E_0$ is the base Young's modulus of the solid material. Typically $p$ is chosen to be $3$ to reduce intermediate density values. Further enhancement of sharp material boundaries is achieved by the addition of a relaxed Heaviside projection \citep{sigmund99new}. In our work, we follow these conventions and set $p=3$; see Sec. \ref{sec:mpm}. It is well known that SIMP supports and can, in practice, generate unrealistic optimal solutions with checkerboard artifacts. A range of smoothing filters have been proposed to alleviate this issue \citep{sigmund88, sigmund2007morphology, sigmund2012sensitivity}; we also apply a density filter and a sensitivity filter in our work, see Sec. \ref{sec:sens}.
\subsection{IPC}

\paragraph{Contact Potential}
\label{sec:contact_potential}

The Incremental Potential Contact (IPC) model \citep{li2020incremental} is a potential-based contact model that guarantees nonpenetration for all configurations. 
For each surface contact pair, the potential is \citep{li2020incremental}
\begin{equation}
    \label{eqn:barrier}
    b(d,\hat{d})=
\begin{cases}
    -\kappa\left(\frac{d}{\hat{d}}-1\right)^2\ln{\left(\frac{d}{\hat{d}}\right)} &0<d<\hat{d} \\
    0 &d\geq \hat{d}
\end{cases},
\end{equation}
where $d$ is the unsigned distance between the two objects in the contact pair (detailed below) and $\hat{d}$ is a user-specified distance threshold, below which the contact potential activates. When $d<\hat{d}$, the barrier energy becomes non-zero and then diverges as $d \to 0,$ allowing it to generate arbitrarily large contact forces, thus preventing penetration ($d\leq0$). Parameter $\kappa$ controls the intensity of the contact force. For a smaller $\kappa$ contact-pair distance $d$ must correspondingly be smaller to generate sufficient contact repulsion. The barrier energy is $C^2$ smooth, ensuring superlinear convergence of Newton's method when solving for displacement $u$. Defining $\mathcal{C}$ as the set of all surface contact pairings, the total contact potential $e_{\text{contact}}$ is then  
\begin{equation}
    \label{eqn:total_contact_potential}
    e_{\text{contact}} = \sum_{k\in \mathcal{C}} \frac{h \hat{d}}{2}b(d_k(x), \hat{d}),
\end{equation}
where $h$ is the discretization's grid spacing, and the weight $\frac{h\hat{d}}{2}$ approximates the integrated energy over the world space \citep{li2022bfemp}. $x$ denotes the world space position of the object. $d_k(x)$ is the distance between the two objects in the contact pair $k.$ Its computation is elaborated below.

\begin{figure}[t]%
\centering
\includegraphics[width=0.42\textwidth]{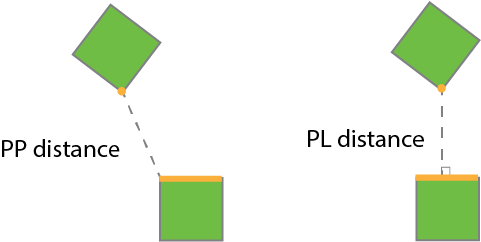}
\caption{Illustrations of point-point (PP) and point-line (PL) distances.}\label{fig:pppe}
\end{figure}

\paragraph{Distances}\label{par:distances} An elastic body is discretized in this work by an axis-aligned regular grid. In 2D, the boundary of a structure consists of vertices and axis-aligned line segments. Each possible contact pair is then a non-incident point-edge pair $\{p,e\}$, where $p$ is a point and $e=(y,z)$ is an edge with endpoints $y$ and $z$. A point $p$ non-incident to edge $e$ implies that both $(p, y)$ and $(p, z)$ do not form boundary edges. The contact distance is then the minimal Euclidean distance between $p$ and $e$, i.e.,
\begin{equation}
\label{eqn:distance}
    d = \min_{\beta}\|p-(y+\beta(z-y))\|\;\; \text{s.t.}\; 0\leq\beta\leq1.
\end{equation}
This constrained optimization problem has two explicit solutions, classified as a point-point (PP) and a point-line (PL) distance. If $\beta=0\text{ or }1$, then the distance is a PP-distance, and 
\begin{equation}
\label{eqn:PP_distance}
    d=\|p-y\|\text{ or }\|p-z\|;
\end{equation}
otherwise, the distance is a PL-distance, and
\begin{equation}
\label{eqn:PE_distance}
    d=\frac{\|(y-p)\times(z-p)\|}{\|y-z\|}.
\end{equation}
These two types of distances are illustrated in Fig.\ref{fig:pppe}.

\paragraph{Friction}
In the IPC framework, friction forces are defined per contact pair. For each contact pair $k,$ a consistently oriented sliding basis $T_{k}(x) \in \mathbb{R}^{d m \times(d-1)}$ is constructed, where $d=2\text{ or }3$ is the problem dimension and $m$ is the number of nodes in the system. The corresponding local frictional force $f_{\text{fric}}^k$ is then defined in terms of   $\tilde{\mathbf{u}}_k$, the local relative sliding displacement orthogonal to the distance gradient, and its corresponding discrete velocity $\tilde{\mathbf{v}}_k = \tilde{\mathbf{u}}_k/\Delta t$. As suggested in \citep{moreau2011unilateral, goyal1991planar, goyal1991planar2}, $f_{\text{fric}}^k$ can be defined by maximizing dissipation rate subject to the Coulomb constraint:
\begin{equation}
f_{\text{fric}}^k(x)=T_{k}(x) \underset{\beta \in \mathbb{R}^{d-1}}{\operatorname{argmin}} \boldsymbol{\beta}^{T} \tilde{\mathbf{v}}_{k} \quad \text { s.t. } \quad\|\boldsymbol{\beta}\| \leq \mu \lambda_{k},
\end{equation}
where $\lambda_k$ is the magnitude of contact force, and $\mu$ is the friction coefficient. This can be equivalently written as 
\begin{equation}
f_{\text{fric}}^k(x)=-\mu \lambda_{k} T_{k}(x) f\left(\left\|\tilde{\mathbf{u}}_{k}\right\|\right) \left(\frac{\tilde{\mathbf{u}}_{k}}{\left\|\tilde{\mathbf{u}}_{k}\right\|}\right),
\end{equation}
where the last term on the right-hand side can be any unit vector if the denominator vanishes. The nonsmooth friction magnitude function $f$ is $1$ if $\left\|\tilde{\mathbf{u}}_{k}\right\|>0$ and falls in $[0,1]$ if the displacement is zero.  \citet{li2020incremental} approximates $f$ by a $C^1$ smooth function 
\begin{equation}
f_{1}(y)=\left\{\begin{array}{ll}
-\frac{y^{2}}{\epsilon_{v}^{2} \Delta t^{2}}+\frac{2 y}{\epsilon_{v} \Delta t}, & y \in\left[0, \Delta t \epsilon_{v}\right) \\
1, & y \geq \Delta t \epsilon_{v}
\end{array}\right.
\label{eq:friction_mollifier}
\end{equation}
where $\epsilon_v$ is a velocity bound such that sliding velocities with magnitude less than $\epsilon_v$ are treated as static.
The friction force is non-integrable. We follow \citet{li2020incremental} to approximate $T$ and $\lambda$ with their values at the previous timestep $t^n$ and $\lambda^n$. Resultingly, the friction force is semi-implicit but can be integrated into a potential energy
\begin{equation}
e_{\text{friction}}^k(x)=\mu \lambda_{k}^{n} f_{0}\left(\left\|u_{k}\right\|\right),
\end{equation}
where $f_0$ satisfies $f_0' = f_1$ and $f_0(\epsilon_v \Delta t) = \epsilon_v \Delta t$ so that $f_{\text{fric}}^k = -\nabla_x e_{\text{friction}}^k.$ The total friction is thus 
\begin{equation}
    f_{\text{fric}} = h^2 \sum_{k\in C} f_{\text{fric}}^k
\end{equation}
and correspondingly, total friction potential can be expressed as 
\begin{equation}
e_{\text{friction}}(x)=h^{2} \sum_{k \in C} e_{\text{friction}}^k(x),
\end{equation}
where $h$ is the mesh spacing and $C$ is the set of all active contact pairs. Here the integration weight per contact pair $k$ is incorporated in $\lambda^k$. See \citep{li2020incremental} for details.

\section{Framework}
\subsection{Material Distribution Representation and Design Variables} \label{sec:mpm}
The design domain $\design$ is discretized with square elements in this work. In SIMP, Young's modulus of each cell $c$ is assumed to be $E_0 \rho^p,$ where $E_0$ is the base Young's modulus of the solid material. Following \citet{rozvany2000simp}, we set $p=3$ to improve binarization of intermediate density values and add a Heaviside projection \citep{sigmund99new}:
\begin{equation}
    H(\rho; \beta) = \frac{\tanh(\frac12\beta) + \tanh(\beta (\rho - \frac12))}{2\tanh(\frac12\beta)}
\end{equation}
to further help with convergence.

Cell densities $\rho_c$ are the design variables for material distribution. We follow \citet{Li_LETO} and solve for force equilibrium via the Material Point Method. Quadrature points within each computational cell share the same density value to avoid subcell QR-pattern artifacts \citep{Li_LETO}. That is, $\rho_q = \rho_c$ for all quadrature points $q$ that belong to the same cell $c.$ 

For a hyperelastic object, the total elastic energy induced by deformation is
\begin{equation}
\label{eqn:elastic_energy}
    e_{\text{elasticity}}(\rho, u) = \int_{\design} \Psi(F)dX,
\end{equation}
where $\Psi$ is the applied energy density function determined by the constitutive model, and $F$ is the deformation gradient,
\begin{equation}
\label{eqn:deformation_gradient}
    F = \frac{\partial x}{\partial X} = I + \frac{\partial u}{\partial X}.
\end{equation}
with $x(X)$ the world space mapping of a material point $X \in \design$ and $u(X)=x(X)-X$ the displacement field.
Following \citet{Li_LETO}, we approximate the elastic energy by 
\begin{equation}
    e_{\text{elasticity}}(\rho, u) \approx E_0 \sum_q  \tilde{\rho}_q^3 \Psi(F_q) V_q,
\end{equation}
where $\tilde{\rho}_q = H(\rho_q; \beta).$
\begin{figure*}[h]%
\centering
\includegraphics[width=0.98\textwidth]{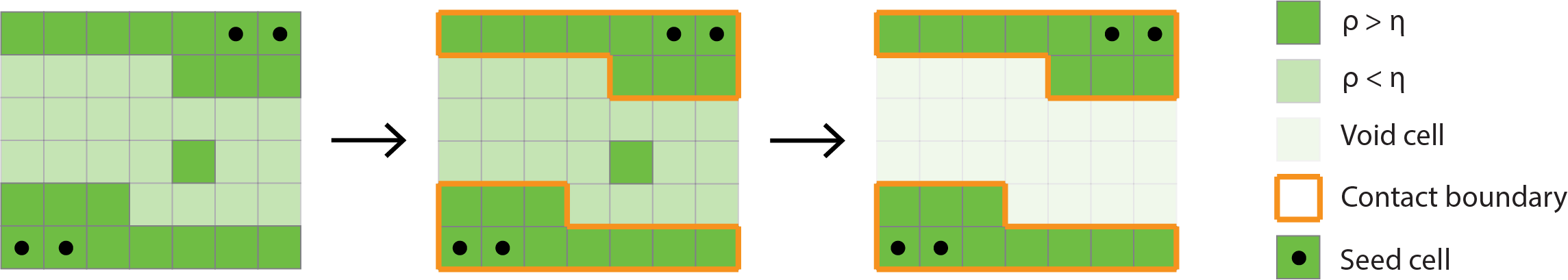}
\caption{The narrow-band procedure. Left: Before applying the narrow-band process. Mid: Identification of largest connected components. Right: Cells outside the established structures and cells with density lower than a threshold $\eta$ are treated as void. The contact boundary is identified when the procedure finishes.}\label{narrowband}
\end{figure*}
Without loss of generality, we apply the compressible neo-Hookean energy density in this work:
\begin{equation}
\label{eqn:NeoHookean_density}
    \Psi(F) = \frac{\mu}{2}\text{tr}(F^T F-d)-\mu\log(J)+\frac{\lambda}{2}\log^2 (J),
\end{equation}
where $J=\det(F)$, $d=2 \text{ or }3$ is the dimension of the problem, and $\mu$ and $\lambda$ are the lam\'e parameters. Note that our numerical procedure guarantees a positive $J$ throughout; see Sec. \ref{sec:pojNewton}.

\subsection{Incorporation of Contact} \label{sec: contact_methodology}
\paragraph{Contact Boundary Detection} \label{para:detect_boundary}

The contact boundary is evolved based on the re-allocated density field during topology optimization iterations, after which it is fed into the contact solver. \citet{bruns2003element} introduced the narrow-band process to topology optimization so that low-density elements can be systematically removed and reintroduced. It has then been widely used, for example, in \citep{liu2018narrow, zhang2021narrow, zhou2016feature}, for filtering out low-density cells and thus avoiding singular stiffness matrices.

Our method adds a boundary detection procedure to the narrow-band process. Let $G$ denote the graph consisting of all cells $c$ defined by the actual adjacency between cells. Given a set of seed cells $S\subseteq G$ (which, by default, is chosen to be where boundary conditions are specified), a depth-first search (DFS) is performed to find all largest connected components $L_i\subseteq G$ such that \begin{equation}
    L_i \cap L_j = \emptyset \text{ if } i\neq j,
\end{equation}
where for each $i,$ \begin{equation}
    \rho_c > \eta, \forall c\in L_i
\end{equation}
and \begin{equation}
    \exists \text{ } c\in S \text{ such that } c \in L_i.
\end{equation}
Here $\eta$ is a thresholding parameter for density. Cells with density lower than $\eta$ and cells that are not connected to any component $L_i$ are removed and treated as zero-density void region. Meanwhile, the boundaries for all components $L_i$ are identified based on mesh connectivity to form the surface boundary for the contact potential. Further, each component $L_i$ is guaranteed to have boundary conditions, thus ensuring a nonsingular stiffness matrix. See Fig. \ref{narrowband} for an illustration of the full procedure.

We remark that minor structures with very small densities can appear in topology optimization iterations that will have no influence on the final optimized material distribution. Here, using the narrow-band process to remove these spurious components can accelerate the convergence of SIMP \citep{liu2018narrow}. Thus, while the narrow-band process enables our auto-detection of the evolving boundary, it also improves the overall convergence of the topology optimization process. 

\paragraph{Contact Potential}
The IPC model is then integrated into the system by utilizing the detected codimension-1 boundary geometries. 
As covered in section \ref{par:distances}, distances are calculated differently for point-point (PP) and point-line (PL) cases. Therefore, we divide the set $\mathcal{C}$ of all the non-incident point-edge pairs into two groups $\mathcal{C}_{PP}$ and $\mathcal{C}_{PL}$ containing only PP pairs and PL pairs, respectively. The total contact potential  (\ref{eqn:total_contact_potential}) is then fully separated as
\begin{align}
\label{eqn:contact_potential}
    e_{\text{contact}}(u) &=\sum_{i\in\mathcal{C}_{PP}}\frac{h\hat{d}}{2}b(d_i,\hat{d})+\sum_{j\in\mathcal{C}_{PL}}\frac{h\hat{d}}{2}b(d_j,\hat{d})
\end{align}
and can be separately evaluated. Classifying the two cases at the energy level ensures the non-ambiguous evaluations of their gradient and Hessian during static solves (See Sec. \ref{sec:pojNewton}). 

\paragraph{Frictional Contact via Artificial Timesteps}

\begin{figure*}
\centering
\begin{subfigure}{.5\textwidth}
  \centering
  \includegraphics[width=.7\linewidth]{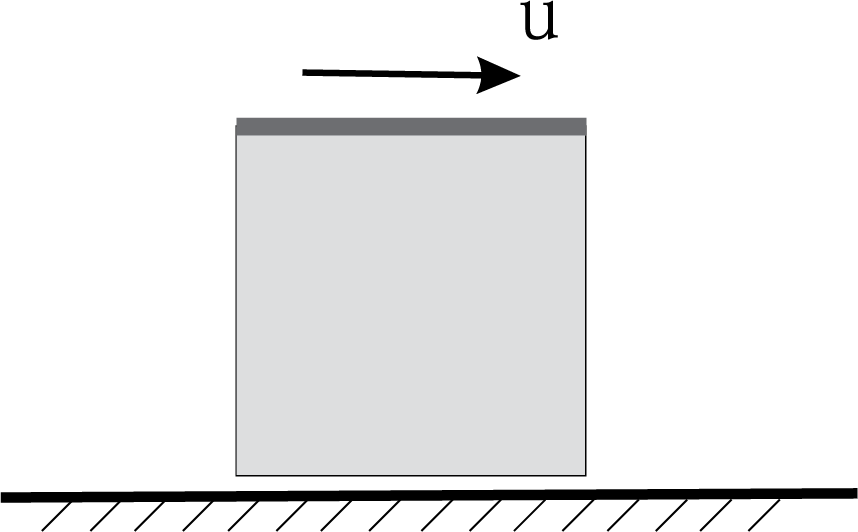}
  \label{fig:sub1}
\end{subfigure}%
\begin{subfigure}{.5\textwidth}
  \centering
  \includegraphics[width=.9\linewidth]{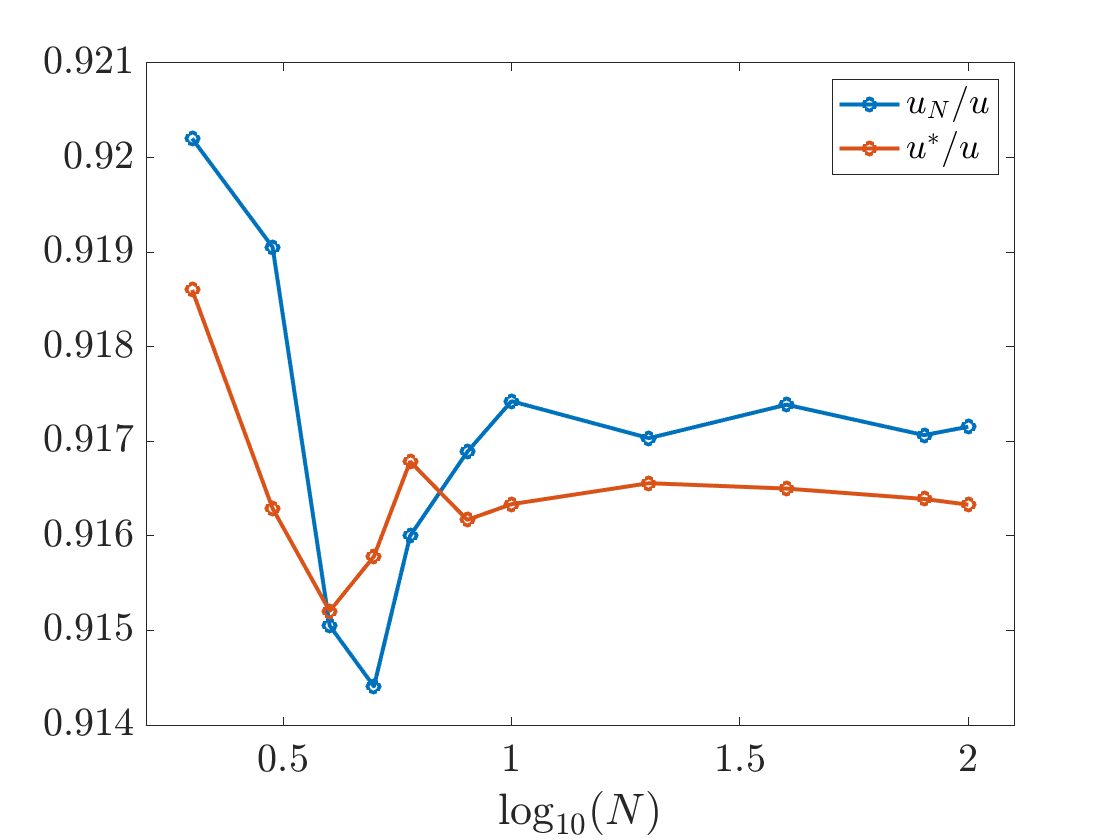}
  \label{fig:sub2}
\end{subfigure}
\caption{An experiment showing the convergence of quasi-static solves under friction. Left: A square of width $20$ cm is aligned against a fixed wall, their distance being $0.1\hat{d}.$ A displacement of $u=20$ cm is prescribed on the top boundary of the square. $T=1$, $\mu = 0.2$ and $\epsilon_v = 10^{-6}.$ Right: A sequence of $N$ quasi-static solves are performed for different values of $N$. Convergence can be observed as $N$ becomes larger. For a balance between accuracy and efficiency, we choose $N=10.$}
\label{fig: friction_experiment}
\end{figure*}

The IPC friction model is parameterized by sliding velocities. While there is no velocity at a static solution, final resting equilibria are supported by friction forces given by the model's sticking conditions. Inspired by \citep{fang2021guaranteed}, we propose to find equilibria under friction via artificial time-stepping. We apply a sequence of $N$ quasi-static solves over an artificial time period $t\in[0,T]$. In the following, we use $T=1$ for all examples. Specifically, we solve a series of nodal displacements $\{u^n\}_{n=1}^N$ such that 
\begingroup\makeatletter\def\f@size{8}\check@mathfonts
\begin{equation}
\begin{split}
    f_{\text{fric}} (u^{n+1}; u^n, \Delta t) + f_{\text{cont}}(u^{n+1}) +f_{\text{int}}(u^{n+1}) + f_{\text{ext}} = 0, \\
    \quad\quad\quad n = 0, 1, ..., N-1
\end{split}
    \label{eqn:friction_equilibrium1}
\end{equation}
\endgroup
 where $u^0=0$ is the displacement at rest, and $\Delta t=\frac{T}{N}$.  Dirichlet boundary displacements are evenly divided into $N$ segments and applied in each corresponding artificial time step. 
 We view $u^n$ and $\Delta t$ as parameters of friction for the definitions of sliding basis and (\ref{eq:friction_mollifier}) respectively. Here $u^N$ (as $N$ is large) can be interpreted as an asymptotic predicted position under dynamic friction via quasi-static approximation. 

 Finally, we solve 
 \begin{equation}
\begin{split}
    f_{\text{fric}} (u; u^{0}, T) + f_{\text{cont}}(x) + f_{\text{int}}(x) + f_{\text{ext}} = 0, 
\end{split}
    \label{eqn:friction_equilibrium}
\end{equation}
 using $u^N$ as the starting point to get $u^*.$ This enables us to find a local minimum $u^*$ satisfying force equilibrium while remaining close to $u^N$. We would like $u^*$ to  converge as $N$ increases so that (\ref{eqn:friction_equilibrium}) has a local minimum independent of $N$. The convergence study of $u^*$ is illustrated in Fig. \ref{fig: friction_experiment}. We also remark that this mechanism only works when solutions to (\ref{eqn:friction_equilibrium}) and every step of (\ref{eqn:friction_equilibrium1}) exist. In our case, we ensure that the Dirichlet boundary condition is defined on each component of the structure.

\subsection{Strain Limiting Relaxation}
\label{sec:strain_limiting}
\begin{figure}[t]%
\centering
\includegraphics[width=0.45\textwidth]{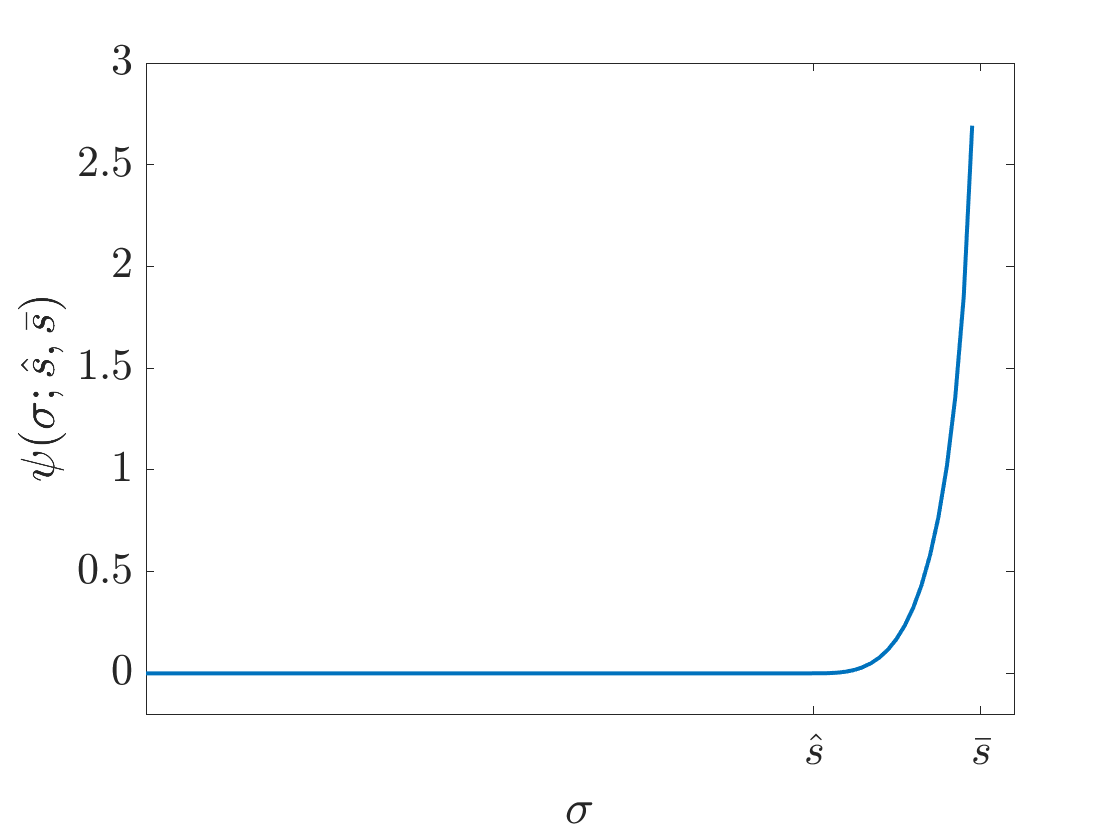}
\caption{$\psi(\sigma; \hat{s}, \bar{s})$ smoothly transits to zero at $\hat{s}$ and grows to infinity at $\bar{s}$.} \label{strain-lmt-psi}
\end{figure}

When large internal contact is present, cells with tiny densities near the contact interface tend to experience extreme deformation. This may cause numerical difficulties for the static solver. To alleviate this issue, we add a strain-limiting energy to those fragile cells to moderate distortion \citep{bridson2002robust, goldenthal2007efficient}. Here we follow \citet{li2020codimensional} and define a $C^2$ scalar function
\begin{equation}
 {\psi}(\sigma; \hat{s}, \bar{s})=\left\{\begin{array}{cc}
-\left(\frac{\hat{s} - \sigma}{\bar{s} - \hat{s}}\right)^{2}  \log \left(\frac{\bar{s}-\sigma}{\bar{s}-\hat{s}}\right), & \sigma \geq \hat{s} \\
0, & \sigma<\hat{s}
\end{array}\right.
\label{strain_lmt_setup}
\end{equation}
so that $\psi(\sigma)>0$ when $\sigma > \hat{s}$ and $\psi(\sigma) \to \infty $ when $\sigma \to \bar{s}$; see Fig. \ref{strain-lmt-psi}.
Given parameters $\hat{s}, \bar{s}, \hat{p}, \text{ and } \bar{p},$ the strain-limiting energy density function ${\Psi}^{\text{SL}}$ is thus defined as 
\begin{align}
    {\Psi}^{\text{SL}}(F) &= {\Psi}^{\text{SL}}(\sigma_i) \\&= \sum_{i = 1}^d \psi(\sigma_i; \hat{s}, \bar{s}) + \psi(-\sigma_;, -\hat{p}, -\bar{p}),
\end{align}
where $\sigma_i$'s are the principal stretches defined by the singular values of $F$. Intuitively, the strain is limited in a way that $\sigma_i$ is not allowed to go beyond $\bar{s}$ or fall below $\bar{p}$. These bounds are guaranteed by the numerical procedure we apply, described in Sec. \ref{sec:pojNewton}.

\begin{figure}[t]%
\centering
\includegraphics[width=0.35\textwidth]{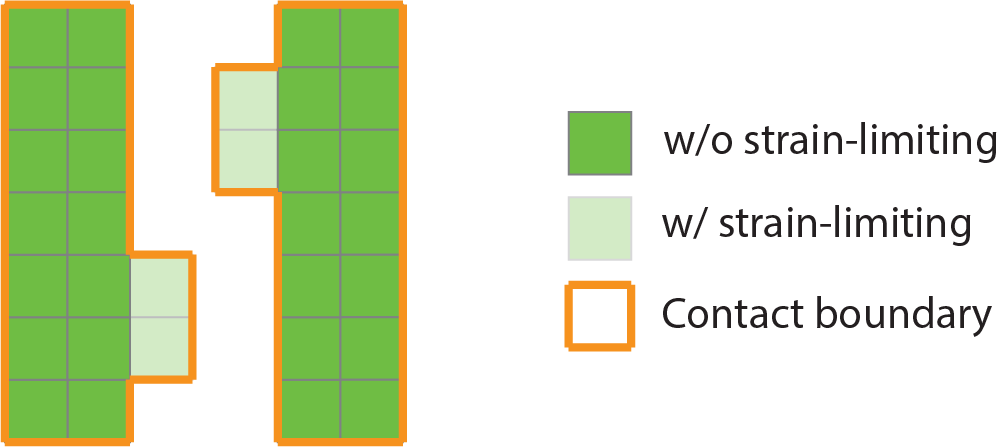}
\caption{Only cells with tiny densities are augmented with strain-limiting relaxation.}\label{fig:strain-lmt}
\end{figure}

Strain limiting is introduced for improving numerical convergence and should not affect the constitutive behavior of well-behaved elements. Therefore, it is only added to cells with densities lower than a specified threshold $\bar{\rho}$; see Fig. \ref{fig:strain-lmt}. Ideally, $\bar{\rho}$ should be as small as possible while numerical convergence is still obtained. We empirically set $\bar{\rho}=0.02.$
To further improve the smoothness of this relaxation, we scale the strain-limiting energy density function with a smooth transition function $h(\cdot)$ such that $h(0) = 1 $ and $h(\bar{\rho}) = 0.$
The total strain-limiting energy can thus be written as \begin{equation}
    e_{\text{SL}} = \sum_{q} h(\rho_q)  \sum_q \Psi^{\text{SL}}(F_q) V_q.
\end{equation}

In practice, we find that a simple linear interpolation \begin{equation}
    h(\rho) = 1 - \frac{\rho}{\bar{\rho}}
\end{equation}
works well.

\begin{figure*}[h]%
\centering
\includegraphics[width=0.95\textwidth]{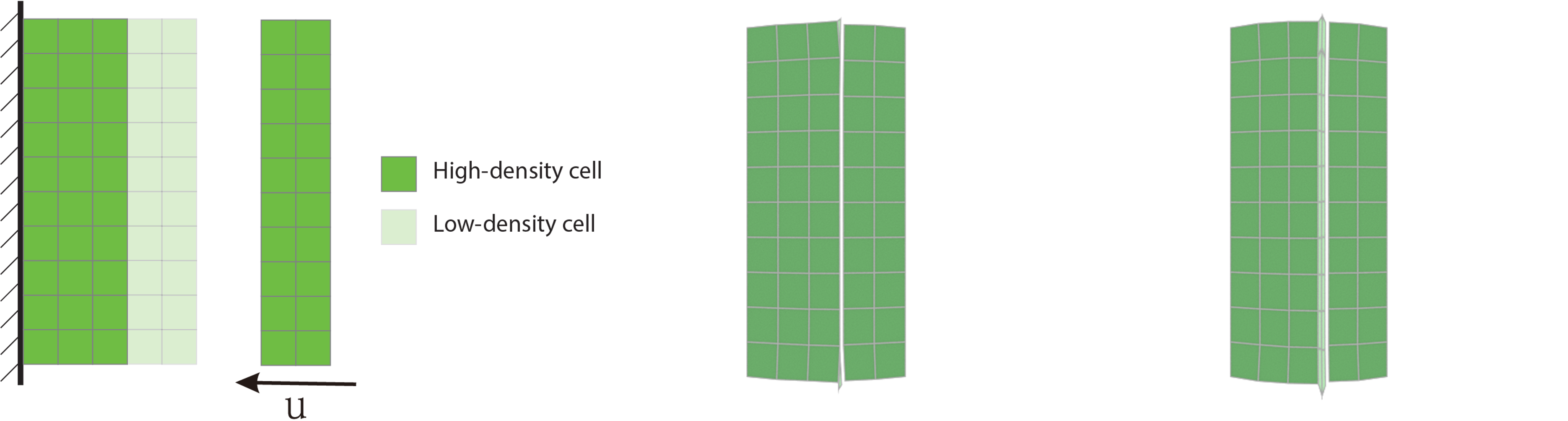}
\caption{An experiment using strain-limiting relaxation. Left: A schematic illustration. Mid: Simulation without strain-limiting. Right: Simulation with strain-limiting. The deformations of cells without relaxation are visually identical in the two cases despite the expected minor quantitative discrepancy.}\label{str_exp}
\end{figure*}

We present in Fig. \ref{str_exp} an experiment to demonstrate that the addition of strain-limiting as relaxation does not significantly alter the deformation of non-relaxed cells.

\subsection{Static Solve with Projected Newton's Method} \label{sec:pojNewton}
Displacement $u$ under static force equilibrium is solved for sensitivity analysis. 

In our framework, as each force (elasticity, the normal contact force, and the tangential frictional force) is associated with a corresponding potential energy, the force equilibrium can be expressed as 
\begin{equation}
    - \frac{\partial e_{\text{total}}}{\partial u} (\rho, u) +  f_{\text{ext}}  = 0, \label{eqn:force-eq}
\end{equation}
where $e_{\text{total}} = e_{\text{elasticity}} + e_{\text{contact}}.$ To incorporate forces due to friction or strain-limiting, it suffices to add the corresponding $e_{\text{friction}}$ or $e_{\text{SL}}$ to $e_{\text{total}}.$
Solving (\ref{eqn:force-eq}) is equivalent to minimizing
\begin{equation}
\label{eqn:min_prob_static}
    \min_{u}e_{\text{total}}(\rho, u)-u \cdot f_{\text{ext}}
\end{equation}
subject to boundary conditions.
Following \citet{Li_LETO}, we use the projected Newton's method to solve the minimization problem (\ref{eqn:min_prob_static}), where the Hessian matrix is projected to be symmetric positive definite (SPD) and a line search procedure is performed to guarantee global convergence \citep{nocedal1999numerical}. Note that for each energy term $e_i$ in $e_{\text{total}},$ a corresponding stepsize upper bound $\alpha_i$ is needed to ensure that the energy $e_i$ is well defined. For instance, in the contact energy $e_{\text{contact}},$ the additive continuous collision detection method \citep{li2020codimensional} is used to bound the stepsize to prevent trajectory intersection. 
For neo-Hookean elasticity and strain-limiting relaxation, an upper bound is derived to prevent the $\log(\cdot)$ term from approaching $\infty$. Finally, the global stepsize upper bound, $\alpha$, is the tightest determined bound, i.e., \begin{equation}
    \alpha = \min_i \{ \alpha_i\}.
\end{equation}

\subsection{Sensitivity Analysis} \label{sec:sens}
Applying the adjoint method, we compute the sensitivity analysis for a general objective function $G(\rho, {u}(\rho))$ with respect to $\rho_c,$ the density of cell $c,$ as 
\begin{equation}
    \frac{dG}{d \rho_c} = \left[\frac{d E_c}{d \rho_c}\right]^T \frac{dG}{d E_c}, 
\end{equation}
 with \begin{equation}
     \frac{dG}{d E_c}  = \frac{\partial G}{ \partial E_c}- \frac{\partial ^2 e}{ \partial E_c \partial {{u}}} \left[\frac{\partial ^2 e }{ \partial {u}^2}\right] ^{-1} \frac{dG}{d{u}},
 \end{equation}
 and \begin{equation}
     \frac{d E_c}{d\rho_c} = 3 E_0 H^2(\rho_c, \beta) \frac{\partial H(\rho_c, \beta)}{\partial \rho_c},
 \end{equation}
 where $E_c$ is defined to be $E_c = E_0 H^3(\rho_c, \beta)$ and $e = e_{\text{total}}$ as in (\ref{eqn:min_prob_static}).
Below we state the sensitivity for three particular objective functions that will be considered in our experiments.\\
\textbf{Compliance} \quad If compliance $G = e_{\text{elasticity}}$ is chosen to be the objective function, then
\begin{equation}
    \frac{\partial G}{\partial E_c} = \sum_{q \text{ in cell } c} V_q \Psi(F_q),
\end{equation}
and $-\frac{\partial G}{\partial {u}}$ is the elastic force.\\

\textbf{Reaction Force} \quad 
The reaction force on node $B$ in direction ${n}$ is defined to be \begin{equation}
    R_{B, {n}} = \frac{\partial e}{\partial u_B} \cdot {n} \in \mathbb{R}.
\end{equation}
Note that the node $B$ must be a Dirichlet node, as a non-Dirichlet node satisfies force equilibrium and hence has zero reaction force. The sensitivity analysis for $R_{B, {n}}$ is therefore 
\begin{equation}
    \frac{d R_{B, {n}}}{d E_c} = \frac{\partial^2 e}{\partial E_c \partial u_B} \cdot {n} -\frac{\partial^2 e}{\partial E_c \partial {u}}\left[\frac{\partial^2 e}{\partial {u}^2}\right]^{-1}\frac{\partial^2 e}{\partial {u} \partial u_B} \cdot {n}.
\end{equation}
Let $D$ denote the set of all nodes where the Dirichlet boundary condition is applied. Forces applied on multiple nodes $\mathscr{B} \subset D$ can be counted together as \begin{equation}
    R = \sum_{B \in \mathscr{B} \subset D} R_{B, {n}}, \quad \frac{dR}{d E_c} = \sum_{B \in \mathscr{B} \subset D} \frac{d R_{B, {n}} }{d E_c}.
\end{equation}
\textbf{Volume Fraction} \quad One important constraint in topology optimization is the volume constraint. The requirement is usually stated as that the total volume fraction does not exceed a threshold $\hat{V}$, i.e., 
\begin{equation}
    g(\rho) = \frac{\sum_c E_c}{\sum_c 1} \leq \hat{V} \in (0,1).
\end{equation}
It follows that \begin{equation}
    \frac{d g}{d E_c} = \frac{E_c}{ \sum_i 1},
\end{equation}
where $i$ loops over all cells. 

\subsection{Density filter, Sensitivity Filter, and Evolving Boundary}
Following \citet{sigmund2007morphology}, density filtering and sensitivity filtering are implemented in our work. Given filter radius $r_{\min},$ the density field $\rho_c$ are modified at the beginning of each iteration:
\begin{equation}
    \widehat{\rho_c} = \frac{1}{\sum_{i \in N_c} H_{ci}} \sum_{i\in N_c} H_{ci}\cdot \rho_i,
\end{equation}
 and the sensitivities $\frac{\partial F}{ \partial \rho}$ are modified before feeding into optimizer:
\begin{equation}
    \widehat{\frac{\partial G}{\partial \rho_c}} = \frac{1}{\max (\gamma, \rho_c) \sum_{i \in N_c} H_{ci}} \sum_{i\in N_c} H_{ci}\cdot \rho_i \frac{\partial G}{\partial \rho_c},
\end{equation}
where \begin{equation}
    H_{ci} = \max ( 0, r_{\min} - \text{dist}(c,i)),
\end{equation}
and \begin{equation}
    N_c = \{i\in \text{all cells}, \text{dist}(c,i) < r_{\min} \}.
\end{equation}
$\gamma = 0.001$ is a small number to avoid division by zero.
\begin{figure*}[h]%
\centering
\includegraphics[width=0.95\textwidth]{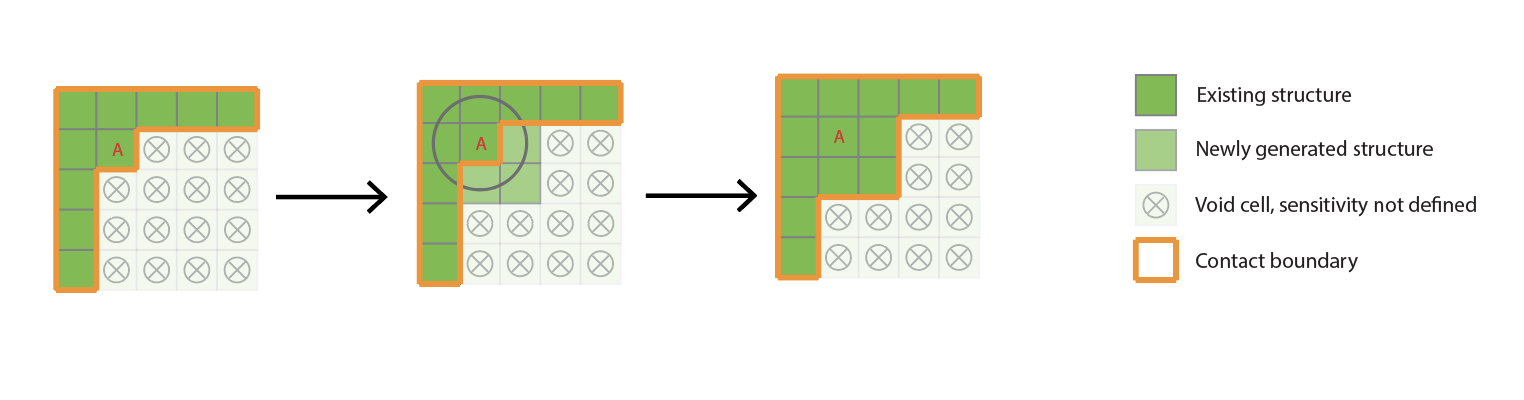} 
\caption{The density and sensitivity filtering and self-evolution of boundary. Sensitivities $dG/d\rho$ are only defined for non-void cells, so they are treated as zero for void cells. Sensitivity filtering and density filtering allow densities to spread from existing cells (e.g. cell $A$) to their neighbours (e.g. cells in circle). The originally void cells are now occupied by the structure, and the contact boundary detected by narrow-band procedure will enlarge correspondingly.}
\label{fig:sens_filter}
\end{figure*}
Intuitively, the filtering functions as a local averaging of densities/sensitivities. As suggested in \citep{sigmund88, sigmund2012sensitivity}, it can ensure the existence of a solution and avoid formations of checkerboard patterns. We follow convention and choose $r_{\text{min}} \in [1.5,3]$ \citep{Li_LETO}. Moreover, in our framework, together with the narrow-band procedure, the filtering allows structures to naturally evolve (both shrink and enlarge) their boundaries, as shown in Fig. \ref{fig:sens_filter}. 

\subsection{Optimizing Structures with MMA} \label{sec:mma}
We adopt the popular optimizer method of moving asymptotes (MMA) \citep{mma_method} to optimize structures. MMA is designed for general structural optimization problems with inequality constraints and box constraints. The original problem is approximated by a series of convex optimizations. At each iteration, two updated asymptotes are set up to constrain the searching interval. MMA typically requires careful parameter tuning for performance. Here we follow the setup in \citep{Li_LETO} and adopt an open-source C++ version of MMA \citep{opensourceMMA}.

\subsection{Overall Pipeline}
The pipeline of this work is summarized in Algorithm \ref{alg:overall}.
\begin{algorithm}[H]
    \caption{Overall pipeline} \label{alg:overall}
    \begin{algorithmic}[1]
    \State {Given: $\rho^0, \hat{V}, f_{\text{ext}}$}
    \For{$j = 0,1,2,3,...$}
        \State Apply narrow-band to detect structure boundaries, see Sec. \ref{para:detect_boundary}
        \State Solve for static equilibrium $u$ via projected Newton, see Algo.\ref{sec:pojNewton}
        \State Compute $\frac{dG}{d\rho}(\rho^j)$, see Sec. \ref{sec:sens}
        \State $\rho^{j+1} \gets \text{MMA}(\rho^j, \frac{dG}{d\rho}(\rho^j))$, see Sec. \ref{sec:mma}
    \EndFor
    \end{algorithmic}
\end{algorithm}

\section{Results}
\subsection{Comparison with Existing Contact Algorithms in Topology Optimization}
\label{sec:other methods}
Here we cover the advantages of the presented method compared to alternatives in topology optimization that handle contact via mortar or fictitious domain methods.

As discussed earlier in Sec. \ref{sec1}, a fundamental drawback of mortar methods is the required pre-specification of contact surfaces prior to simulation and so optimization. This restricts its application only to simple examples with just external contact obstacles where contact surfaces can be easily predetermined and assumed to be unchanging. Thus, once specified, contact boundaries are then generally unable to evolve, disallowing any large changes to be accounted for. Further difficulties then also arise when we require modeling of multiple contacting domains and/or when the accurate and precise resolution of contact requires intersection-free geometries. It is thus challenging to apply mortar methods to handle intricate and often changing internal contacts along evolving boundaries \citep{stromberg2013influence, kristiansen2020topology}.

\begin{figure}
        \centering
        \begin{subfigure}[b]{0.23\textwidth}
            \centering
            \includegraphics[width=\textwidth]{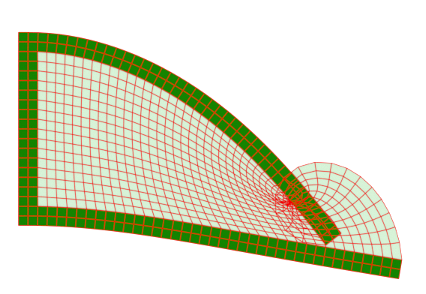}
            \caption[]%
            {{\small }}    
            \label{fig:fig_squeeze_airmesh}
        \end{subfigure}
        \hfill
        \begin{subfigure}[b]{0.23\textwidth}  
            \centering 
            \includegraphics[width=\textwidth]{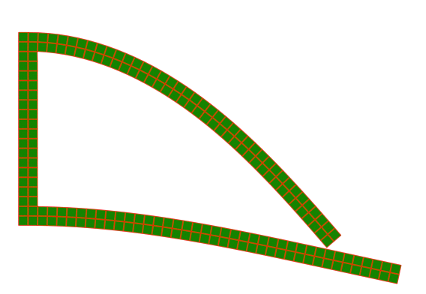}
            \caption[]%
            {{\small }}    
            \label{fig:fig_squeeze_ipc}
        \end{subfigure}
        \vskip\baselineskip
        \begin{subfigure}[b]{0.23\textwidth}   
            \centering 
            \includegraphics[width=\textwidth]{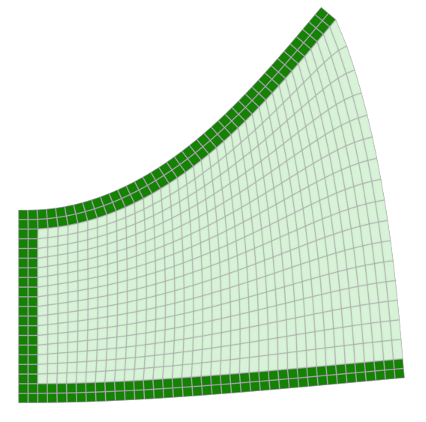}
            \caption[]%
            {{\small}}    
            \label{fig:fig_drag_airmesh}
        \end{subfigure}
        \hfill
        \begin{subfigure}[b]{0.23\textwidth}   
            \centering 
            \includegraphics[width=\textwidth]{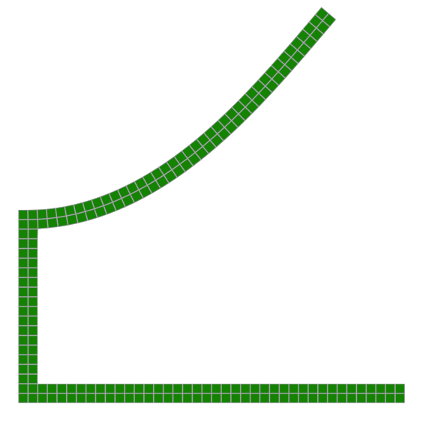}
            \caption[]%
            {{\small}}    
            \label{fig:fig_drag_ipc}
        \end{subfigure}
        \caption[]
        {\small (a) and (b): top branch is dragged downwards, simulated with air-mesh and IPC, respectively. (c) and (d): top branch is pulled upwards, simulated with air-mesh and IPC, respectively. Observe that in (c), the bottom branch is slightly elevated, a phenomenon that should not occur.} 
        \label{fig:mean and std of nets}
    \end{figure}

Fictitious domain methods, on the other hand, are prone to numerical artifacts. As a comparison, we simulate a C-shaped structure deformed by pulling its top right corner downwards. Fig.s \ref{fig:fig_squeeze_airmesh} and \ref{fig:fig_squeeze_ipc} demonstrate the converged simulation results for this set-up using respectively an air-mesh and IPC model. 
The air-mesh model requires per-example fine-tuning of stiffness parameters to ensure that the resulting gap between geometries that should be in contact is neither too large (and so not really in contact yet) nor negative (to avoid self-penetration artifacts) \citep{Wriggers_third_medium, weissenfels2015contact}. Moreover, the inversion of the domain and/or significant shearing of the fictitious domain is generally unavoidable. \citet{sigmund_internal_contact} can partially alleviate difficulties for air-mesh models by wrapping the domain with further layers of air-mesh, but these challenges remain.

Convergence of the forward simulation and accurate satisfaction of non-intersecting geometries are guaranteed independent of the choice of $\hat{d}$. In turn, the parameter $\hat{d}$ gives direct control of how close materials can be prior to application of contact forces (see dashed lines in Fig. \ref{airmesh_v_ipc_scheme}). This enables users to decide how accurately conforming contact geometries should be per application. As smaller $\hat{d}$ increases accuracy at the cost of more computation, it provides a parameter for directly controlling accuracy versus efficiency. In a similar manner, the convergence parameter for the Newton solve itself then allows a choice to balance efficiency versus accuracy for the simulation solves.   

Another disadvantage of the air-mesh is that the fictitious domain elements generate non-physical artificial forces (and thus inaccurate deformations) when deformations of the domain are large. For example, when two boundaries are pulled sufficiently far away from each other, we see this effect in even simple examples like Fig. \ref{fig:fig_drag_airmesh}. Clearly, these regions without material (voids) should not apply forces on the domain in these contexts. On the other hand, the IPC model, e.g. as simulated in Fig. \ref{fig:fig_drag_ipc}, ensures that forces are solely applied between true material surfaces and so avoids these artifacts altogether. In addition, by extending the simulated region and so increasing its range of deformation, we observe that air-mesh models also incur additional computational challenges. In our experiments, we see that for a single static solve, the air-mesh method generally takes about five to ten times more Newton iterations when compared to the IPC method to converge to the same tolerance for the same example. Finally, when it comes to modeling frictional contact air-mesh models lack the appropriate resolution of the necessary terms to model tangential resistance under contact. Friction modeling thus remains a major challenge for fictitious domain methods. Here, in contrast, the IPC model includes a direct and natural friction model (Sec. \ref{sec: contact_methodology} and Sec. \ref{sec:friction_result}) that ensures proper coupling between well-defined normal forces and tangential friction forces.

\subsection{Fixed-Interface Contact}
\begin{figure}[h]%
\centering
\includegraphics[width=0.45\textwidth]{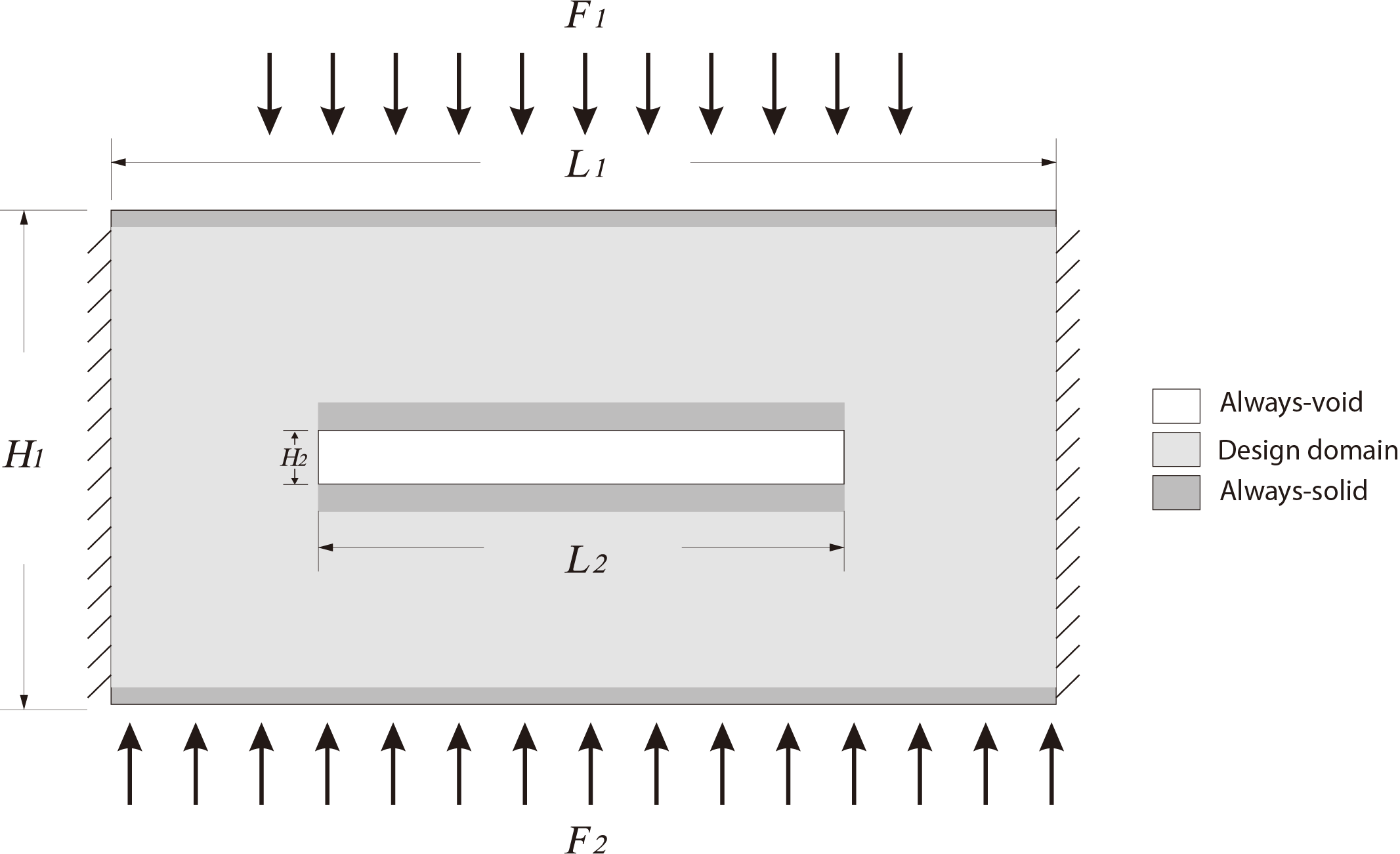}
\caption{Problem setup for fixed-interface contact problem.}\label{dd2}
\end{figure}

\begin{figure}[h]
     \centering
     \begin{subfigure}[b]{0.45\textwidth}
         \centering
         \includegraphics[width=\textwidth]{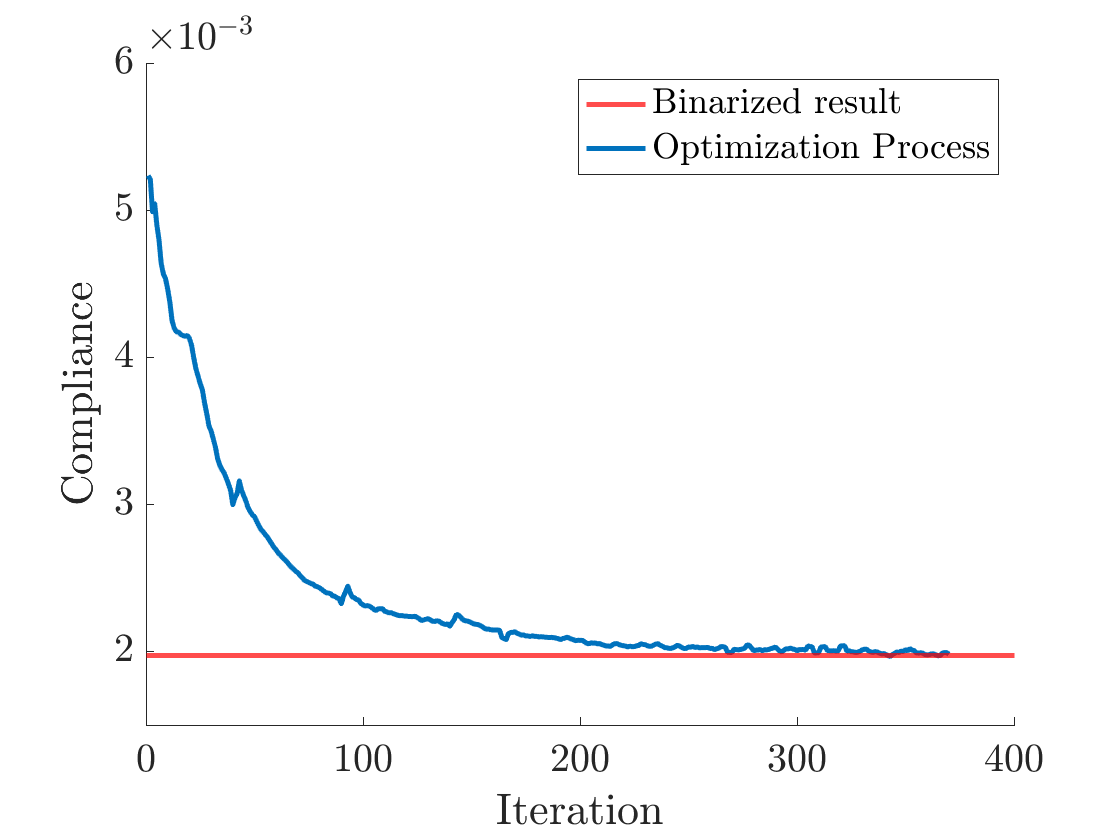}
         \caption{}
         \label{fig:midgap_compliance}
     \end{subfigure}
     \hfill
     \begin{subfigure}[b]{0.45\textwidth}
         \centering
         \includegraphics[width=\textwidth]{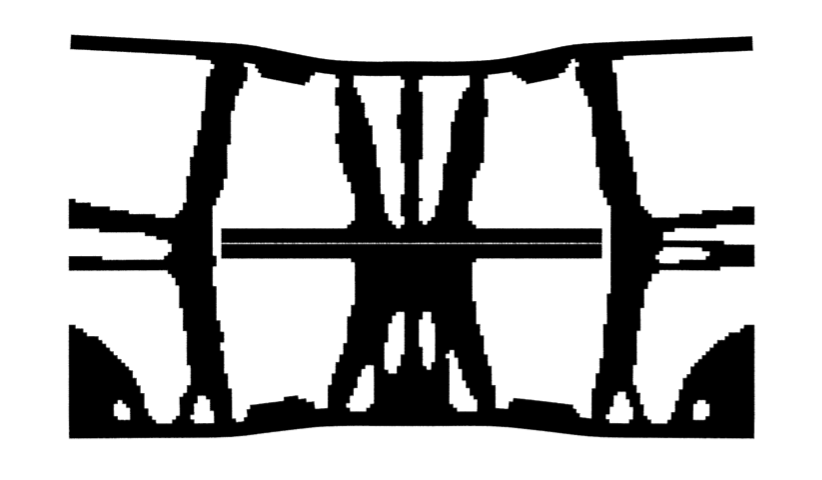}
         \caption{}
         \label{fig:midgap_disp}
     \end{subfigure}
     \caption[]
        {\small(a): Quantitative compliance plot. The compliance of the binarized system is marked in red. (b): World space displacement of the binarized result.} 
\end{figure}

We first test our method in an experiment designed to ensure internal contacts to the design domain will occur at pre-specified, fixed interface. Importantly, as our results will demonstrate, a resulting optimal design can still consider and include contact interfaces in other regions as well. For the experimental set up, please see Fig. \ref{dd2}. 
Here the design domain is of width $L_1 = 1.8m$ and height $H_1 = 1m.$ An internal ``always-void" region of width $L_2 = 1m$ and height $H_2 = 0.08m$ is applied in the center where we keep $\rho = 0$ fixed throughout. Downward forces of $F_1 = 0.72N$ are evenly loaded on the top surface, length $\frac{2}{3}L_1$, and corresponding upward forces of $F_2 = 1N$ are evenly loaded on the entire bottom face surface of the domain. A homogeneous Dirichlet boundary condition is enforced along both the two side walls. Simulation resolution is $180 \times 80,$ with a mesh resolution of $h = 0.01m.$ The IPC parameter $\hat{d}$ is $0.1h,$ $\mu = 0$ and $E_0$ is set to $100$. To facilitate the formation of contact supporting structures, two additional ``always-solid'' regions (marked in dark grey) are specified above and beneath the ``always-void" region with a fixed density of $\rho = 1$. We also set another two thin layers of ``always-solid" regions: one at the top and one at the bottom of the domain. Non-fixed portions of the design domain are then initialized with $\rho = 0.42$ with a  volume constraint of $42\%.$ Here, we begin with a goal objective of minimizing structure compliance with equilibrium given by the force balance of elastic and contact forces.
\begin{figure}
     \centering
     \begin{subfigure}[b]{0.39\textwidth}
         \centering
         \includegraphics[width=\textwidth]{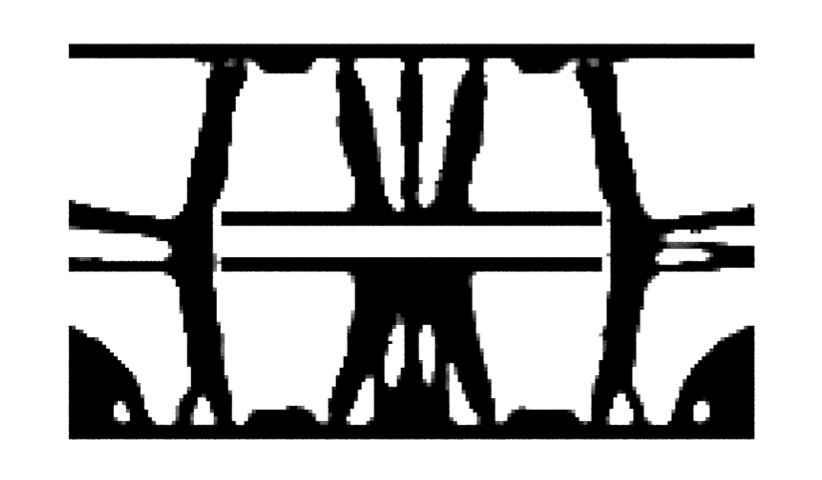}
         \caption{Density field after 370 iterations.}
         \label{fig:midgap_final}
     \end{subfigure}
     \hfill
     \begin{subfigure}[b]{0.39\textwidth}
         \centering
         \includegraphics[width=\textwidth]{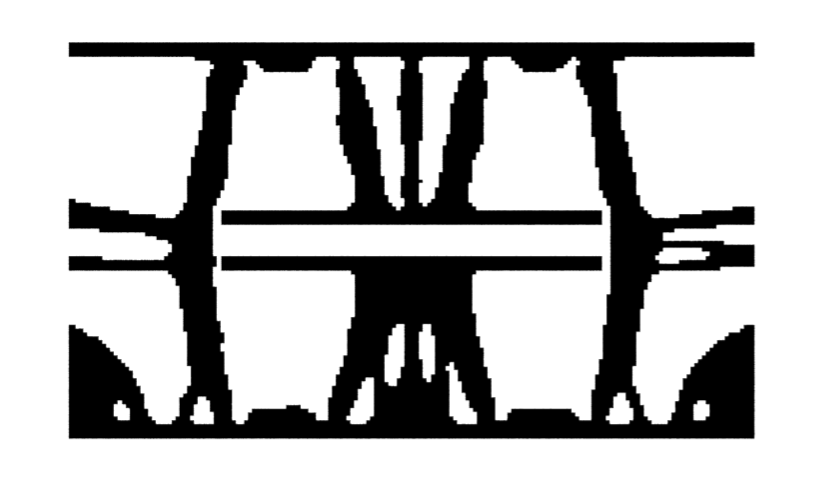}
         \caption{Density field after binarization.}
         \label{fig:midgap_polar}
     \end{subfigure}
     \caption[]
        {\small Optimized result for the fixed-interface contact problem.} 
\end{figure}
Convergence can be observed after around 370 iterations, where both the density field (Fig. \ref{fig:midgap_final}) and the value of the objective function (Fig. \ref{fig:midgap_compliance}) approach invariance. 

A final binarization \begin{equation}
    B(\tilde{\rho}) = \left\{\begin{array}{cc}
1 \text{ if } \tilde{\rho}  \geq 0.5\\
0 \text{ if } \tilde{\rho } < 0.5
\end{array}\right.
\end{equation}
is then applied to get the binary 0-1 solution shown in Fig. \ref{fig:midgap_polar}.
Compliance of the optimized solution decreases from a start of $5.234 \times 10^{-3}$ to $1.981\times 10^{-3},$ while the binarized solution reduces a bit further to $1.976 \times 10^{-3}$ (the red horizontal line in Fig. \ref{fig:midgap_compliance}). The volume constraint $V(\rho) / \hat{V} \times 100\%$ is $100.09\%$ for the converged result and $100.01\%$ for the binarized result. 
It can be observed in Fig. \ref{fig:midgap_disp} that two surfaces designed for contact here do indeed touch closely along their entire interface. Here the final optimized deformed structure is primarily supported by two vertical beams connecting the top plate and the bottom plate and by contacts along the ``fixed" contact interface. Peripheral structures then also connect the left and right (where the homogeneous Dirichlet condition is specified) to the bottom plate and to the vertical beams for additional support.

\subsection{Two-Stage Min-Max Problem}
\begin{figure}[h]%
\centering
\includegraphics[width=0.45\textwidth]{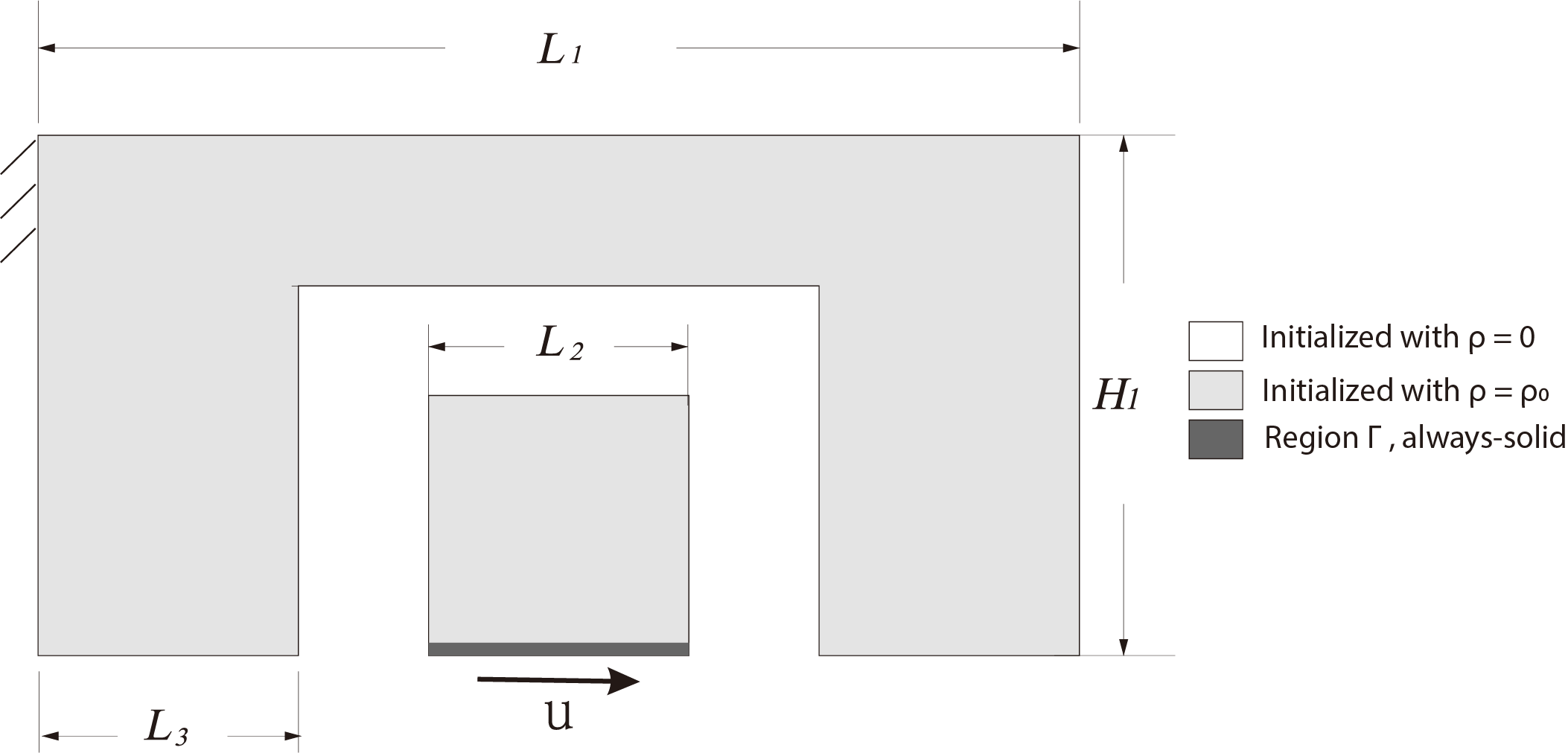}
\caption{Design domain of the two-stage min-max problem.}\label{dd1}
\end{figure}
Next, we apply our contact-aware topology optimization algorithm to design a structure that will handle switching from a loose configuration to close contact under varying magnitudes of a prescribed displacement. The problem setup is shown in Fig. \ref{dd1}. Here the design domain has length $L_1 = 1.0m$ and height $H_1 = 0.5m$, with an inset square of width $L_3 = \frac{2}{15}L_1$ and and supporting legs with bottom length $L_2 = \frac{6}{15}L_1.$ The gap, initially set to be void, separates the inner square and the outer piece. The bottom layer of the inner square $\Gamma$ is set to be ``always-solid", and a displacement ${u}$ is prescribed for it. Here the design goal for this problem is to minimize the reaction force on $\Gamma$ for small displacement ${u}$ and, at the same time, to maximize reaction forces when ${u}$ is large. Specifically, the applied design objective is \begin{equation}
    \min_{\rho} \left(R_{\Gamma, {n}} (\rho, {u} = u_1 {n})\right)^2 - \left(R_{\Gamma, {n}} (\rho, {u} = u_2 {n})\right)^2,
\end{equation}
where ${n} = (1,0)$ is the unit vector pointing to the right with small and large displacements respectively $u_1 = \frac{1}{9}L_1,$ and $u_2 = \frac{2}{9}L_1.$ Here the grey region in Fig. \ref{dd1} is initialized with $\rho = 0.32,$ the applied volume constraint is $32\%,$, with $E_0=100, \mu = 0$, and a discretization resolution of $180\times 90.$ Our problem setup thus motivates from a comparable design experiment in \citep{sigmund_internal_contact}. Severe internal contact induced by the larger (but not smaller) displacement of the self-evolving boundary highlights the difficulty of this problem. Along with elasticity and the contact forces, we also apply the above-described strain-limiting relaxation to help reduce numerical difficulty from the large deformations. 
\begin{figure}[h]
        \centering
        \begin{subfigure}[b]{0.23\textwidth}
            \centering
            \includegraphics[width=\textwidth]{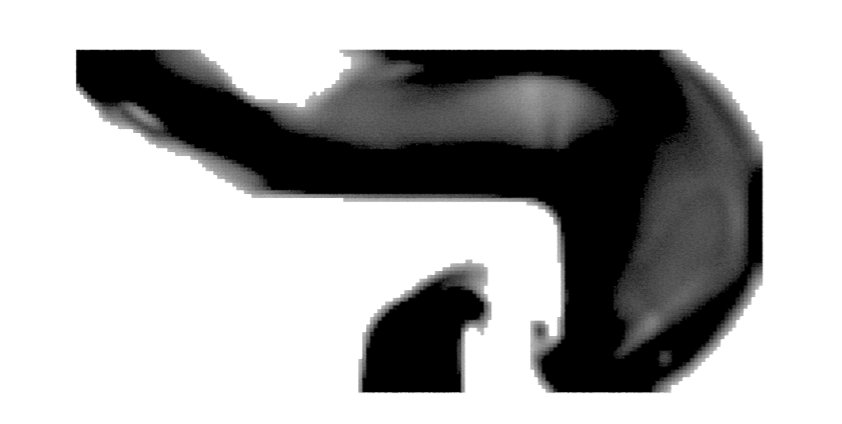}
            \caption[]%
            {{\small After 290 iterations}}    
        \end{subfigure}
        \hfill
        \begin{subfigure}[b]{0.23\textwidth}  
            \centering 
            \includegraphics[width=\textwidth]{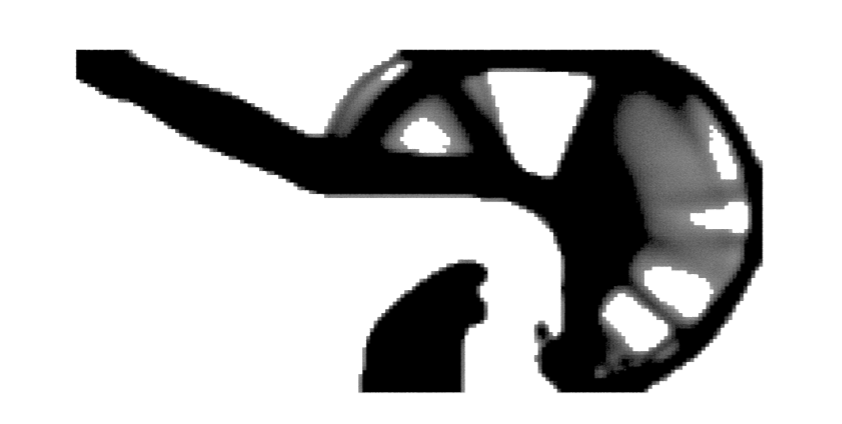}
            \caption[]%
            {{\small After 580 iterations}}    
        \end{subfigure}
        \vskip\baselineskip
        \begin{subfigure}[b]{0.23\textwidth}   
            \centering 
            \includegraphics[width=\textwidth]{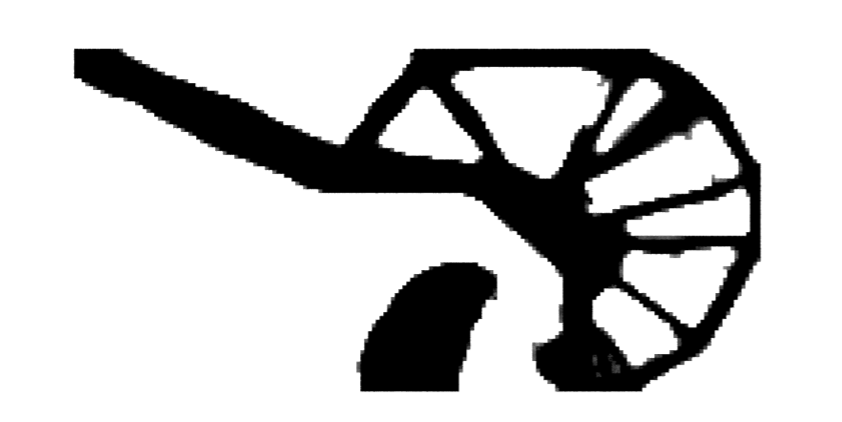}
            \caption[]%
            {{\small After 850 iterations}}   
        \end{subfigure}
        \hfill
        \begin{subfigure}[b]{0.23\textwidth}   
            \centering 
            \includegraphics[width=\textwidth]{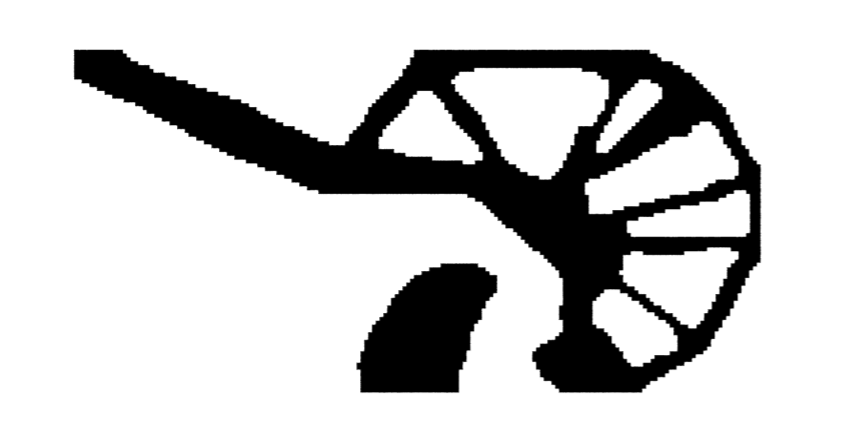}
            \caption[]%
            {{\small Binarized result}}    
            \label{fig:sigmund_polarized}
        \end{subfigure}
        \caption[]
        {\small Density field of the two-stage min-max problem at different stages.} 
        \label{fig:sigmund_stages}
    \end{figure}
\backmatter
Convergence can be obtained after $850$ iterations. We plot the density field after 290 iterations, 580 iterations, and 850 iterations, together with the binarized structure in Fig. \ref{fig:sigmund_stages} to show the optimization process. 

\begin{figure}[h]
     \centering
     \begin{subfigure}[b]{0.45\textwidth}
         \centering
         \includegraphics[width=\textwidth]{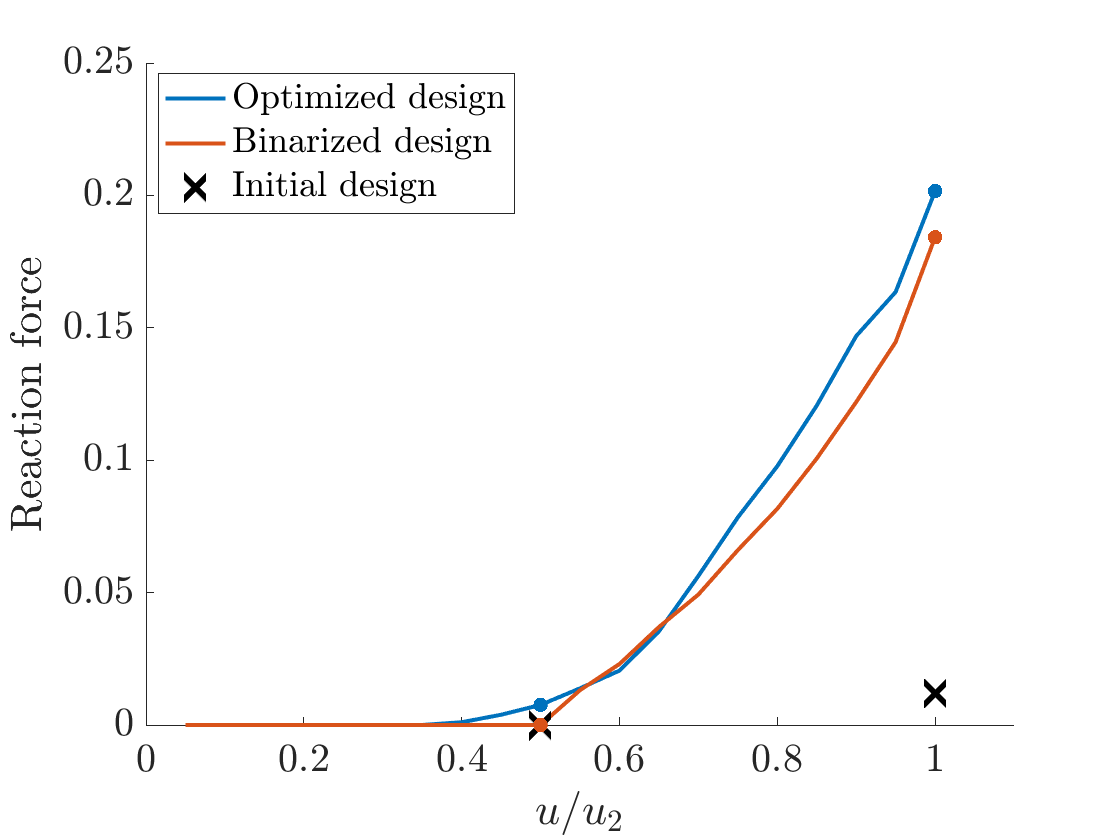}
         \caption{}
         \label{fig:sig_obj}
     \end{subfigure}
     \hfill
     \begin{subfigure}[b]{0.45\textwidth}
         \centering
         \includegraphics[width=\textwidth]{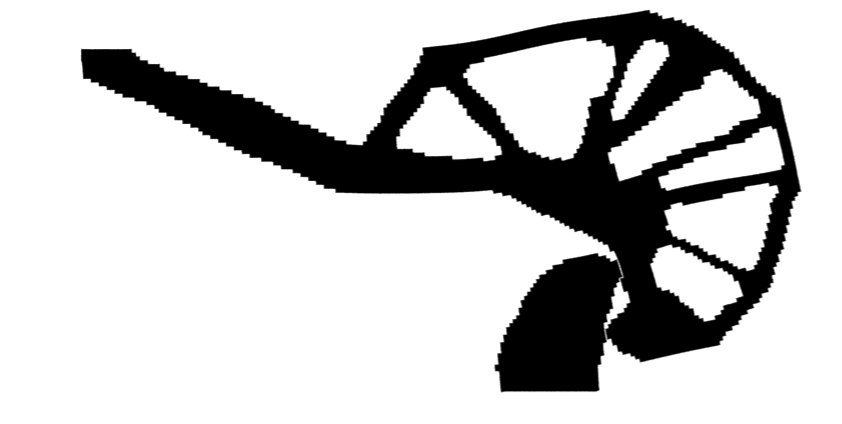}
         \caption{}
         \label{fig:sig_disp}
     \end{subfigure}
     \caption[]
        {\small(a): Reaction force on $\Gamma$ for different level of displacement. $u/u_2 = 0.5$ corresponds to scenario where reaction force is minimized, and $u/u_2 = 1.0$ corresponds to where reaction force is maximized. (b): World space displacement at $u=u_2$ of the binarized result.} 
\end{figure}

Notice from Fig. \ref{fig:sigmund_polarized} that the inner structure and the outer structure, as expected, are separated by a small gap. This reflects the trade-off that the optimizer takes when balancing reaction forces under different values of displacement. The reaction force is minimized when a smaller displacement $u_1$ is prescribed. Ideally the two pieces should remain fairly isolated under prescribed displacement $u_1$ so that there will be no deformation and hence zero reaction force. On the other hand, the reaction force is maximized when a larger displacement $u_2$ is prescribed. Ideally the two pieces should be completely connected so that the prescribed displacement can yield the largest deformation (hence the largest reaction force). Balancing the two goals, the optimal design shown in Fig. \ref{fig:sigmund_polarized} separates the inner piece and the outer piece by a gap just enough to achieve zero reaction force under displacement $u_1,$ while allowing contact to take place for $u>u_1.$ See the red curve in Fig. \ref{fig:sig_obj}. Also, note that compared with the initial design (reaction forces under the two levels of displacement are marked by black crosses), both the optimized design and the final binarized design have achieved significant gain in optimizing the objective.

Fig. \ref{fig:sig_obj} plots the reaction force on $\Gamma$ as the displacement specified on $\Gamma$ ranges from zero to $u_2 \Vec{n}$ for the optimized design and the binarized design. The reaction forces corresponding to $u_1$ and $u_2$ are accentuated, while the reaction forces for the initial design at $u_1$ and $u_2$ are also marked. For the initial design, there is no contact at $u_1,$ resulting in a zero reaction force. At $u_2$ the reaction force is merely $0.0120N.$ The optimized design has significant gain at $u_2$ where the reaction force becomes $0.2016N,$ with a little sacrifice at $u_1$ where there is now slight contact. After removing the peripheral low-density cells, the binarized result re-achieves the state of no contact at $u_1.$ Reaction force at $u_2$ is subsequently reduced to $0.1842N,$ though still much larger than that in the initial design. $V(\rho) / V$ is $100.21\%$ for the optimized design and $100.92\%$ for the binarized design.

\subsection{Screwdriver with Friction} \label{sec:friction_result}
\begin{figure}[h]%
\centering
\includegraphics[width=0.47\textwidth]{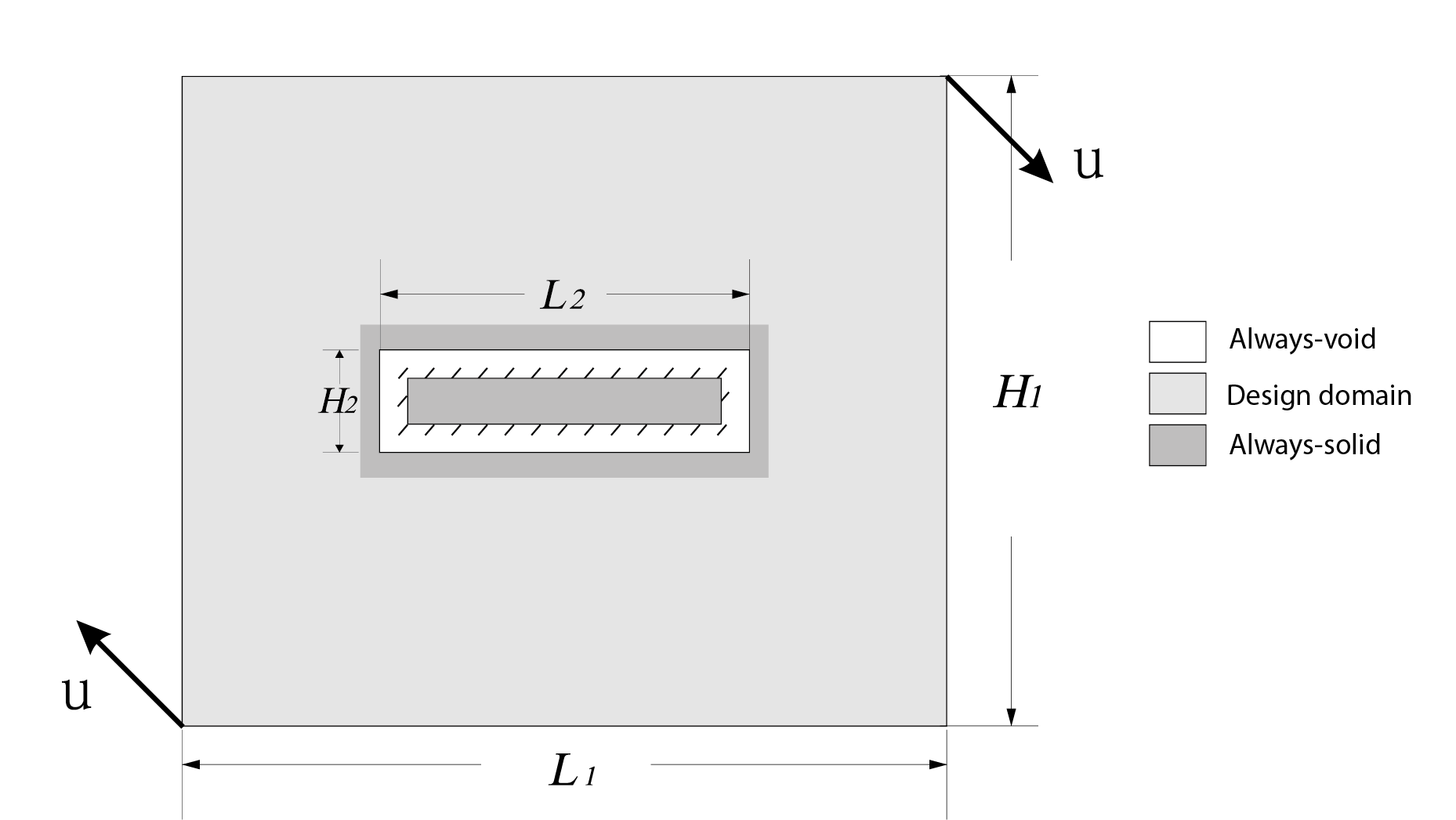}
\caption{Problem setup for the screwdriver problem.}\label{dd3}
\end{figure}
Lastly, we explore the effect of friction using the proposed algorithm. The problem setup is shown in Fig. \ref{dd3}. The design domain is the outer geometry, representing a screwdriver that is rotated to drive the screw inside. The inner rectangle, representing a screw, is fixed in position and remains an ``always-solid" region throughout the optimization process. A gap of width $1.5\hat{d}$ is set between the inner and the outer geometries. (Recall that two objects are treated as in contact when their gap is less than $\hat{d}.$) In this experiment, we set the design domain with $L_1 = 92cm, H_1 = 78cm, L_2 = 44cm,$ and $H_2 = 10cm.$ and a simulation resolution of $92 \times 78.$ Here, the inner layer of the design domain is also kept solid throughout the optimization to ensure the design obtains sliding contact between the screw and the screwdriver. The design domain (light grey region in Fig. \ref{dd3}) is initialized with $\rho = 0.32$ with a volume constraint is set to be $32\%.$ The rotation is prescribed by displacements of $(-16cm, 16cm)$ at the lower left corner and $(16cm, -16cm)$ the upper right corner while $E_0$ is set to 100. Here our design goal objective is to maximize the compliance of the final structure.

For optimizing, we find that incrementally allocating the prescribed displacement ${u}$ into ten static solves (so each progressing displacement by $\frac{1}{10} {u}$) is sufficient to obtain convergent displacement field with friction. See Sec. \ref{sec: contact_methodology} for details.

\newcolumntype{C}{>{\centering\arraybackslash}m{11em}}
\begin{table*}[t]
\centering
\begin{tabular}{l*3{C}@{}} 
      \toprule
             & $\mu=0$ & $\mu=0.2$ & $\mu=0.4$\\
            \midrule
            \small{Optimized design} & \includegraphics[width=11em]{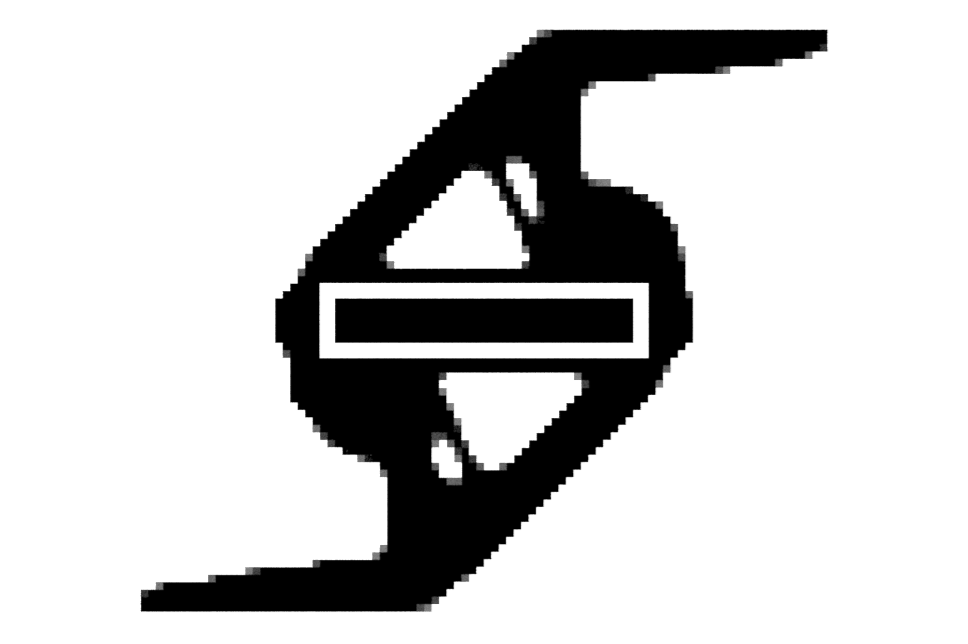} & \includegraphics[width=11em]{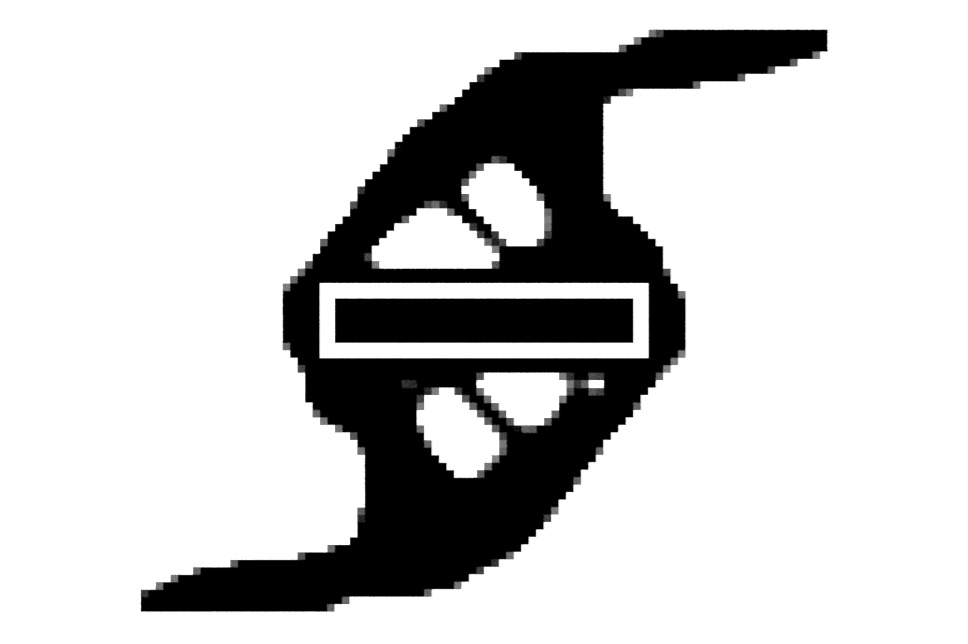}  & \includegraphics[width=11em]{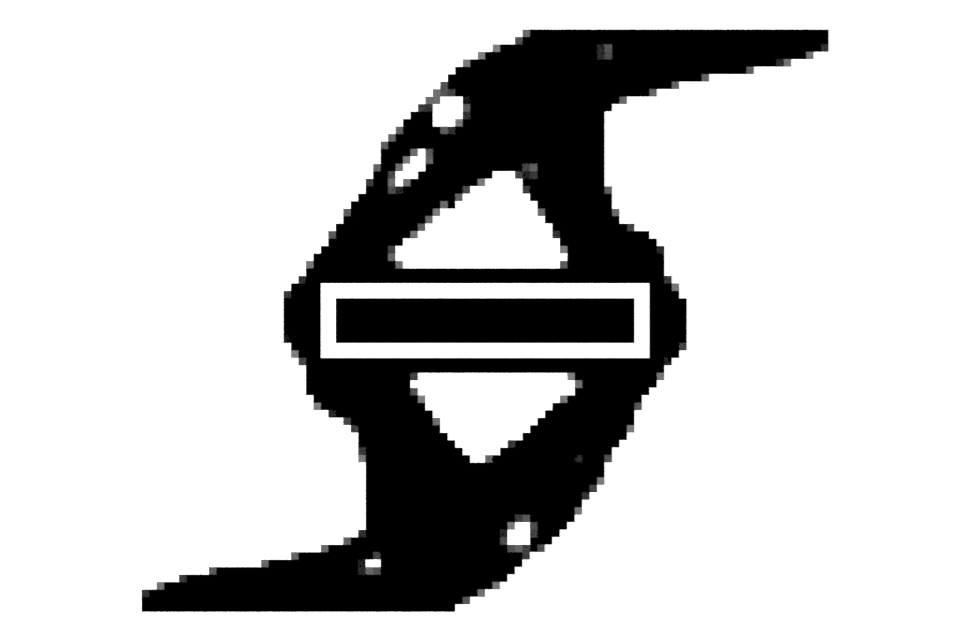}  \\
            \small{Binarized design} & \includegraphics[width=11em]{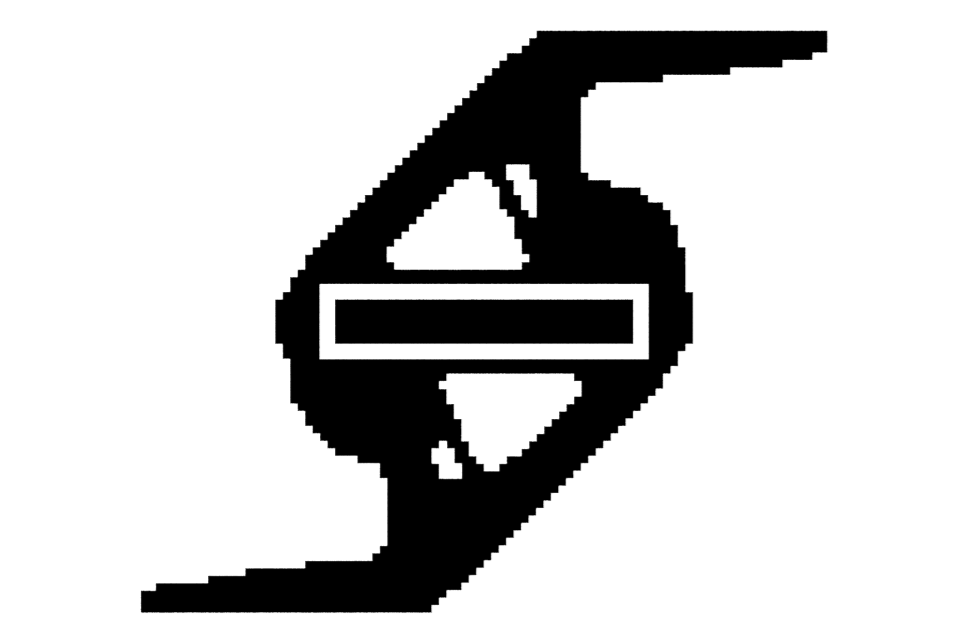}  & \includegraphics[width=11em]{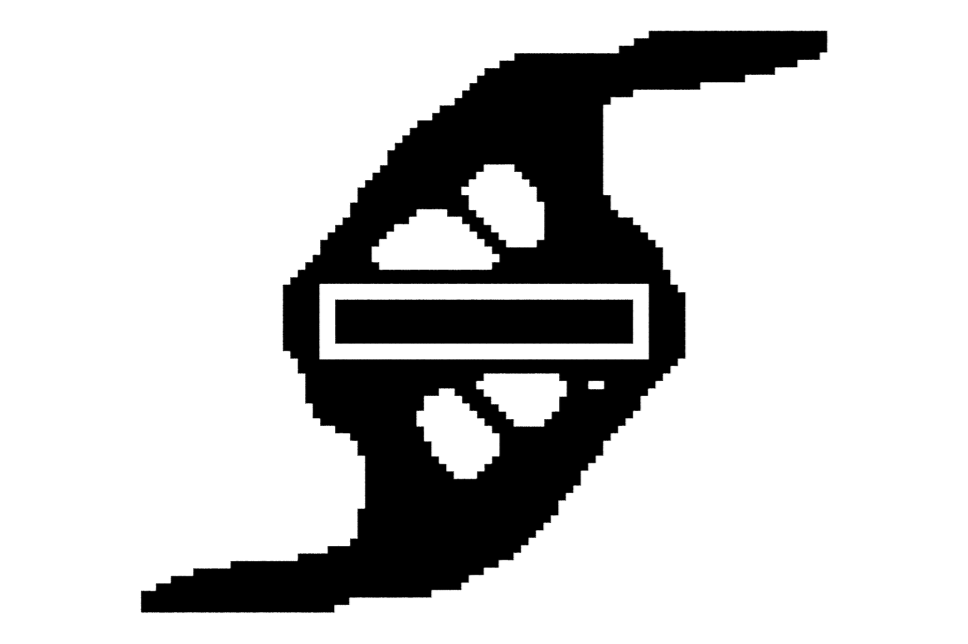}  & \includegraphics[width=11em]{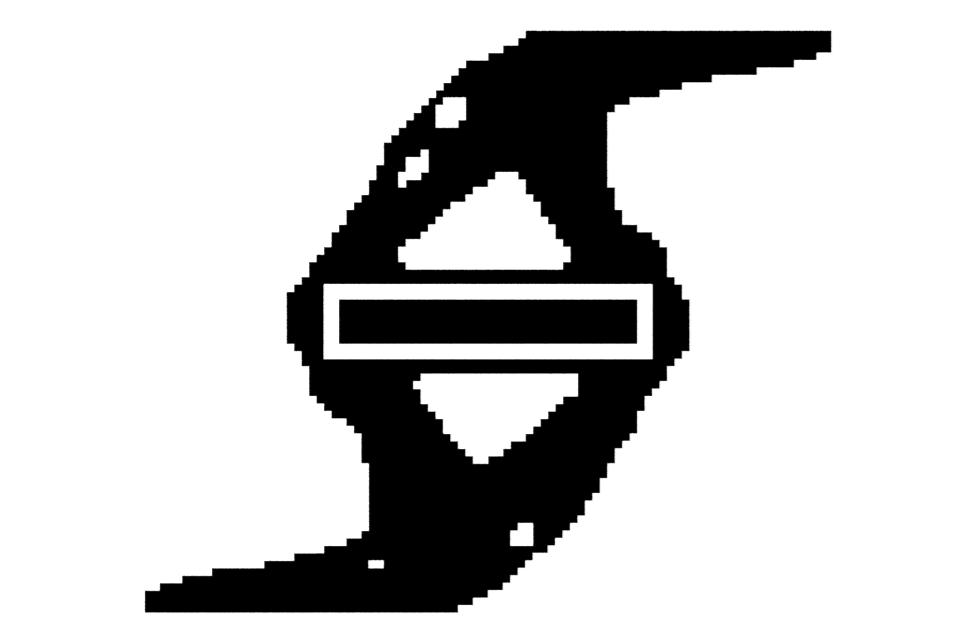} \\
            \small{Deformed} & \includegraphics[width=11em]{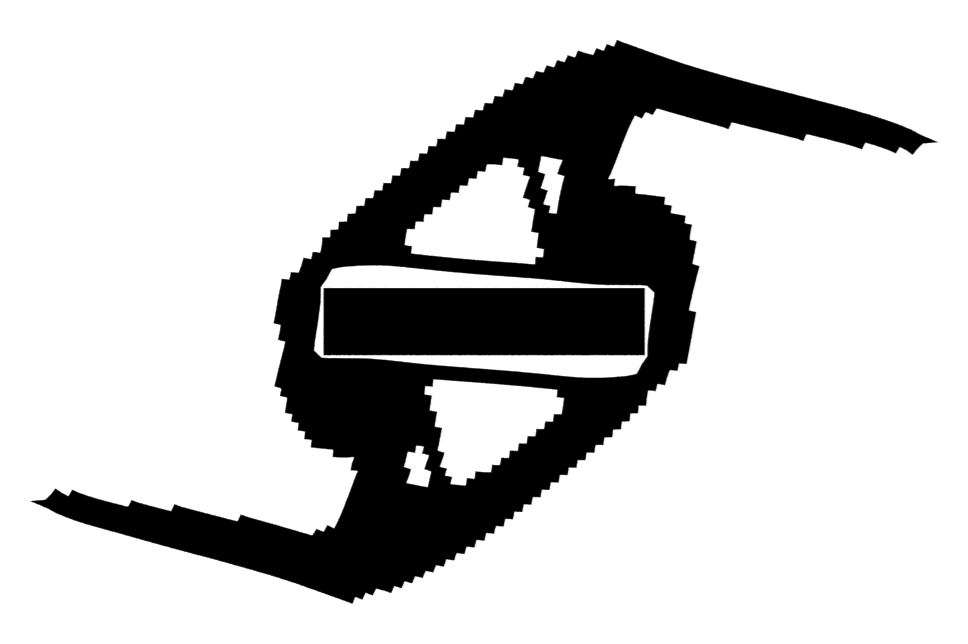}  & \includegraphics[width=11em]{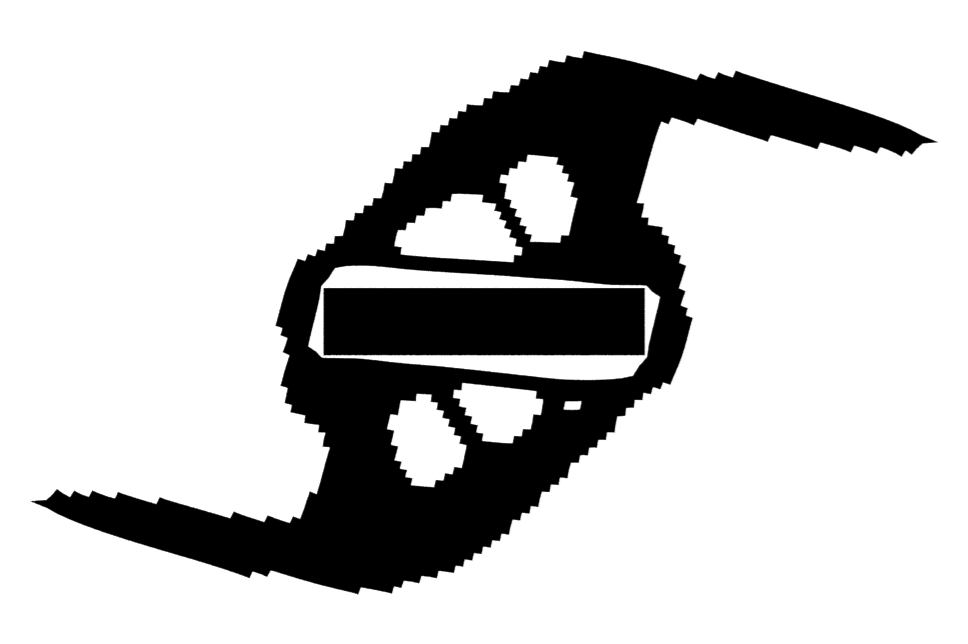}  & \includegraphics[width=11em]{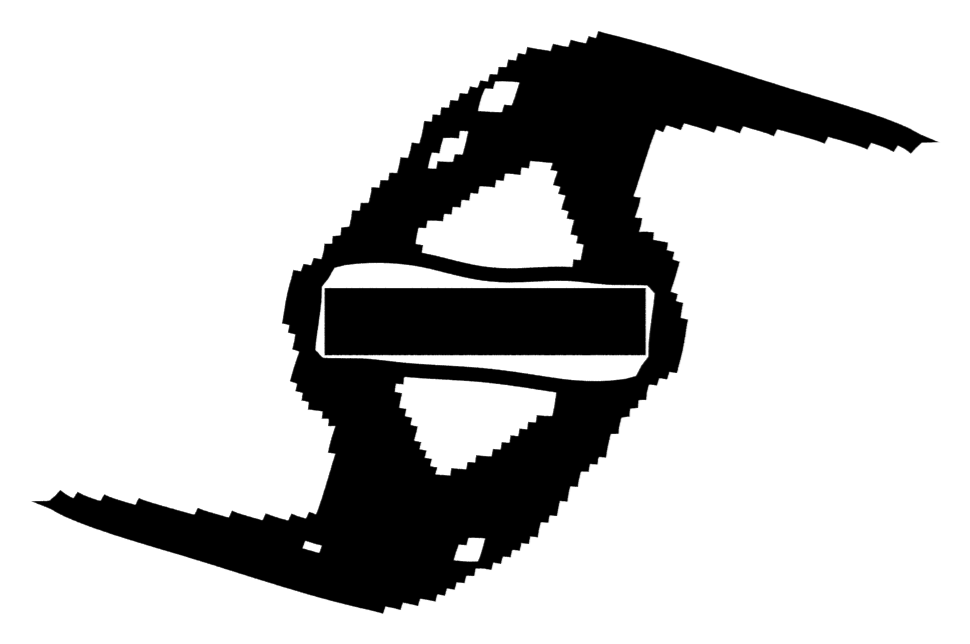}   \\
            \bottomrule 
        \end{tabular}
\caption{Results for different values of friction coefficient $\mu.$ The top row shows the optimized design. The mid row shows the binarized design. The bottom row shows the world space deformation of binarized design.} \label{fig:friction_result}
\end{table*} 
    
\begin{figure}[h]%
\centering
\includegraphics[width=0.45\textwidth]{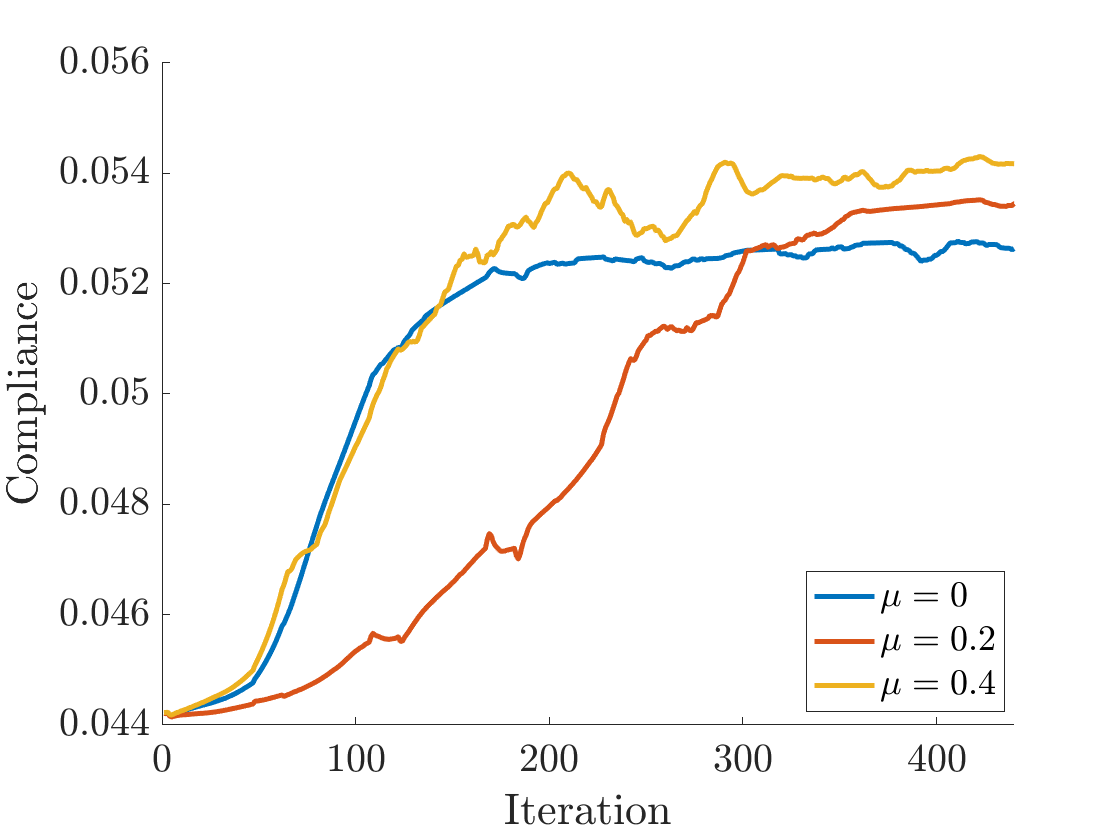}
\caption{Compliance of the structure for different friction coefficients.}\label{fig:friction_obj}
\end{figure}

We solve the design optimization problem for three increasing values of friction: $\mu = 0, \mu = 0.2,$ and $\mu = 0.4.$ 
$440$ iterations were run, and convergence can be observed in all three scenarios. Results are summarized in Table \ref{fig:friction_result}, with change in  compliance demonstrated in Fig. \ref{fig:friction_obj}. Post-evaluation reveals that the compliance of a binarized structure (first row in Fig. \ref{fig:friction_result}) differs no more than 2\% from the compliance of a corresponding optimized result (second row in Table \ref{fig:friction_result}). The volume constraints for the optimized results are reported to be $99.98\%, 99.98\%,$ and $100.77\%$ for $\mu = 0, \mu = 0.2,$ and $\mu = 0.4.$ The volume constraints for the binarized results are reported to be $99.96\%, 100.01\%,$ and $100.73\%,$ respectively. Both the objective values and the volume constraints confirm that convergence has been reached in each case.

Several observations can be made from these results. First, in the last row of Table \ref{fig:friction_result} we see that structures corresponding to larger friction coefficients $\mu$ have less relative displacement at the contact interface. This demonstrates that the friction model, and so downstream optimization, captures the effect of larger $\mu$ applying larger impedance. Second, we see that larger friction values generate optimized designs with more and smaller poles, while with no friction, we see even fewer and larger poles. One possible explanation for this effect is that compared with the case of no friction, over-clustering of mass around where contact happens will yield larger contact forces and hence larger local frictional impedance. This, in turn, will limit deformation to smaller regions. On the other hand, a more uniform distribution of supporting structures will allow deformation to be spread out and can hence achieve a larger compliance. Finally, as we see in Table \ref{fig:friction_result}, a larger friction coefficient $\mu$ generates an optimal design with larger compliance. This is because when the same boundary condition is applied on the lower left and upper right corners, larger impedance due to larger friction yields more distortion within the ``screwdriver,'' and hence larger compliance. 

\section{Discussion}
We have proposed a new framework to handle frictional self-contact in topology optimization. To do so, the IPC model is incorporated into the SIMP algorithm. 
The presented method provides the first frictional contact-aware topology optimization framework with guaranteed non-interpenetration satisfaction covering both external- and self-contact without the need for pre-specification or labeling. 
As demonstrated, this framework now enables optimizations to explore evolving contact interfaces, and so new designs of structures can be generated that are able to take advantage of self-contact. Potential applications of this framework include soft robotic grippers, energy-absorbing cushions, and meta-materials with microstructures in contact.

While the narrow-band procedure was originally introduced to topology optimization to speed up convergence, here, the proposed method augments it with contact boundary detection. This procedure, nevertheless, is then non-smooth in $\rho.$ When compared with a smooth boundary, the grid-aligned contact boundary generated by the narrow-band procedure may also lead to an unnatural concentration of contact in small regions. Thus a promising direction for future work is to form a differentiable and so smooth contact surface using techniques such as \citep{meshSDF}. Last but not least, our current framework uses MMA for design optimization, which usually requires parameter tuning for high performance. In our experiments, a default standard MMA setup is applied without further fine-tuning. Exploring the effects of alternate optimization methods as well as a thorough analysis of their parameters on the current method should be a useful direction for further improvements in practical performance.

 \bmhead{Acknowledgments}
The work was supported in part by the National Science
Foundation of the United States under funding numbers 2011471, 2016414, 2153851, 2153863, 2023780.

\section*{Declarations}

\bmhead{Conflict of interest} The authors declare that there are not competing interests or conflict of interests. 

\bmhead{Ethical approval} This project does not contain any studies with human participants or animals.

\bmhead{Replication of results}
The presented methodology is implemented in C++. The program is compiled with the GNU C++ compiler and executed on the Ubuntu OS. The code and data are freely available to readers upon request.

\bibliography{1_reference}

\begin{thebibliography}{41}
\providecommand{\natexlab}[1]{#1}
\providecommand{\url}[1]{{#1}}
\providecommand{\urlprefix}{URL }
\providecommand{\doi}[1]{\url{https://doi.org/#1}}
\providecommand{\eprint}[2][]{\url{#2}}
 \bibcommenthead

\bibitem[{Andreassen et~al(2011)Andreassen, Clausen, Schevenels, Lazarov, and
  Sigmund}]{sigmund88}
Andreassen E, Clausen A, Schevenels M, et~al (2011) Efficient topology
  optimization in matlab using 88 lines of code. Structural and
  Multidisciplinary Optimization 43(1):1--16

\bibitem[{Bluhm et~al(2021)Bluhm, Sigmund, and
  Poulios}]{sigmund_internal_contact}
Bluhm GL, Sigmund O, Poulios K (2021) Internal contact modeling for finite
  strain topology optimization. Computational Mechanics 67:1099–--1114

\bibitem[{Bridson et~al(2002)Bridson, Fedkiw, and Anderson}]{bridson2002robust}
Bridson R, Fedkiw R, Anderson J (2002) Robust treatment of collisions, contact
  and friction for cloth animation. In: Proceedings of the 29th annual
  conference on Computer graphics and interactive techniques, pp 594--603

\bibitem[{Bruns and Tortorelli(2003)}]{bruns2003element}
Bruns TE, Tortorelli DA (2003) An element removal and reintroduction strategy
  for the topology optimization of structures and compliant mechanisms.
  International journal for numerical methods in engineering 57(10):1413--1430

\bibitem[{Dumas(2018)}]{opensourceMMA}
Dumas J (2018) Mma and gcmma. \url{https://github.com/jdumas/mma}

\bibitem[{Fang et~al(2021)Fang, Li, Jiang, and Kaufman}]{fang2021guaranteed}
Fang Y, Li M, Jiang C, et~al (2021) Guaranteed globally injective 3d
  deformation processing. ACM Trans Graph 40(4):75--1

\bibitem[{Fernandez et~al(2020)Fernandez, Puso, Solberg, and
  Tortorelli}]{Fernandez_mortar_2020}
Fernandez F, Puso MA, Solberg J, et~al (2020) Topology optimization of multiple
  deformable bodies in contact with large deformations. Computer Methods in
  Applied Mechanics and Engineering 371:113,288

\bibitem[{Ferrari and Sigmund(2020)}]{sigmund99new}
Ferrari F, Sigmund O (2020) A new generation 99 line matlab code for compliance
  topology optimization and its extension to 3d. Structural and
  Multidisciplinary Optimization 62(4):2211--2228

\bibitem[{Goldenthal et~al(2007)Goldenthal, Harmon, Fattal, Bercovier, and
  Grinspun}]{goldenthal2007efficient}
Goldenthal R, Harmon D, Fattal R, et~al (2007) Efficient simulation of
  inextensible cloth. In: ACM SIGGRAPH 2007 papers. p 49--es

\bibitem[{Goyal et~al(1991{\natexlab{a}})Goyal, Ruina, and
  Papadopoulos}]{goyal1991planar}
Goyal S, Ruina A, Papadopoulos J (1991{\natexlab{a}}) Planar sliding with dry
  friction part 1. limit surface and moment function. Wear 143(2):307--330

\bibitem[{Goyal et~al(1991{\natexlab{b}})Goyal, Ruina, and
  Papadopoulos}]{goyal1991planar2}
Goyal S, Ruina A, Papadopoulos J (1991{\natexlab{b}}) Planar sliding with dry
  friction part 2. dynamics of motion. Wear 143(2):331--352

\bibitem[{Hallquist et~al(1985)Hallquist, Goudreau, and
  Benson}]{hallquist1985sliding}
Hallquist J, Goudreau G, Benson D (1985) Sliding interfaces with contact-impact
  in large-scale lagrangian computations. Computer methods in applied mechanics
  and engineering 51(1-3):107--137

\bibitem[{Han et~al(2022)Han, Xu, Duan, and Huang}]{Han_mortar}
Han Y, Xu B, Duan Z, et~al (2022) Stress‐based topology optimization of
  continuum structures for the elastic contact problems with friction.
  Structural and Multidisciplinary Optimization 65(2)

\bibitem[{Kristiansen et~al(2020)Kristiansen, Poulios, and
  Aage}]{kristiansen2020topology}
Kristiansen H, Poulios K, Aage N (2020) Topology optimization for compliance
  and contact pressure distribution in structural problems with friction.
  Computer Methods in Applied Mechanics and Engineering 364:112,915

\bibitem[{Kruse et~al(2018)Kruse, Nguyen-Thanh, Wriggers, and
  De~Lorenzis}]{kruse2018isogeometric}
Kruse R, Nguyen-Thanh N, Wriggers P, et~al (2018) Isogeometric frictionless
  contact analysis with the third medium method. Computational Mechanics
  62(5):1009--1021

\bibitem[{Li et~al(2020{\natexlab{a}})Li, Ferguson, Schneider, Langlois, Zorin,
  Panozzo, Jiang, and Kaufman}]{li2020incremental}
Li M, Ferguson Z, Schneider T, et~al (2020{\natexlab{a}}) Incremental potential
  contact: intersection-and inversion-free, large-deformation dynamics. ACM
  Trans Graph 39(4):49

\bibitem[{Li et~al(2020{\natexlab{b}})Li, Kaufman, and
  Jiang}]{li2020codimensional}
Li M, Kaufman DM, Jiang C (2020{\natexlab{b}}) Codimensional incremental
  potential contact. arXiv preprint arXiv:201204457

\bibitem[{Li et~al(2021{\natexlab{a}})Li, McWilliams, Li, Sung, and
  Jiang}]{SoftHAV}
Li X, McWilliams J, Li M, et~al (2021{\natexlab{a}}) Soft hybrid aerial vehicle
  via bistable mechanism. In: 2021 IEEE International Conference on Robotics
  and Automation (ICRA), pp 7107--7113, \doi{10.1109/ICRA48506.2021.9561434}

\bibitem[{Li et~al(2022)Li, Fang, Li, and Jiang}]{li2022bfemp}
Li X, Fang Y, Li M, et~al (2022) Bfemp: Interpenetration-free mpm--fem coupling
  with barrier contact. Computer Methods in Applied Mechanics and Engineering
  390:114,350

\bibitem[{Li et~al(2021{\natexlab{b}})Li, Li, Li, Zhu, Zhu, and
  Jiang}]{Li_LETO}
Li Y, Li X, Li M, et~al (2021{\natexlab{b}}) Lagrangian--eulerian multidensity
  topology optimization with the material point method. International Journal
  for Numerical Methods in Engineering 122(14):3400--3424

\bibitem[{Liu et~al(2018)Liu, Hu, Zhu, Matusik, and Sifakis}]{liu2018narrow}
Liu H, Hu Y, Zhu B, et~al (2018) Narrow-band topology optimization on a
  sparsely populated grid. ACM Transactions on Graphics (TOG) 37(6):1--14

\bibitem[{Luo et~al(2016)Luo, Li, Duan, and Kang}]{Luo_mortar_hyperelastic}
Luo Y, Li M, Duan Z, et~al (2016) Topology optimization of hyperelastic
  structures with frictionless contact supports. International Journal of
  Solids and Structures 81:373--382. \doi{10.1016/j.ijsolstr.2015.12.018}

\bibitem[{Mankame and Ananthasuresh(2004)}]{mankame2004topology}
Mankame ND, Ananthasuresh G (2004) Topology optimization for synthesis of
  contact-aided compliant mechanisms using regularized contact modeling.
  Computers \& structures 82(15-16):1267--1290

\bibitem[{Moreau(2011)}]{moreau2011unilateral}
Moreau JJ (2011) On unilateral constraints, friction and plasticity. In: New
  variational techniques in mathematical physics. Springer, p 171--322

\bibitem[{Müller et~al(2015)Müller, Chentanez, Kim, and
  Macklin}]{Muller_fictitious}
Müller M, Chentanez N, Kim T, et~al (2015) Air meshes for robust collision
  handling. ACM Transactions on Graphics (TOG) 34(4):1--9

\bibitem[{Niu et~al(2019)Niu, Zhang, and Gao}]{niu2019topology}
Niu C, Zhang W, Gao T (2019) Topology optimization of continuum structures for
  the uniformity of contact pressures. Structural and Multidisciplinary
  Optimization 60(1):185--210

\bibitem[{Nocedal and Wright(1999)}]{nocedal1999numerical}
Nocedal J, Wright SJ (1999) Numerical optimization. Springer

\bibitem[{Pagano and Alart(2008)}]{Pagano_fictitious}
Pagano S, Alart P (2008) Self-contact and fictitious domain using a difference
  convex approach. International journal for numerical methods in engineering
  75(1):29--42

\bibitem[{Remelli et~al(2020)Remelli, Lukoianov, Richter, Guillard,
  Bagautdinov, Baqu{\'{e}}, and Fua}]{meshSDF}
Remelli E, Lukoianov A, Richter SR, et~al (2020) Meshsdf: Differentiable
  iso-surface extraction. CoRR abs/2006.03997.
  {\href{https://arxiv.org/abs/2006.03997}{{https://arxiv.org/abs/2006.03997}}}

\bibitem[{Rozvany(2000)}]{rozvany2000simp}
Rozvany G (2000) The simp method in topology optimization-theoretical
  background, advantages and new applications. In: 8th Symposium on
  Multidisciplinary Analysis and Optimization, p 4738

\bibitem[{Sigmund(1997)}]{sigmund1997design}
Sigmund O (1997) On the design of compliant mechanisms using topology
  optimization. Journal of Structural Mechanics 25(4):493--524

\bibitem[{Sigmund(2001)}]{sigmund200199}
Sigmund O (2001) A 99 line topology optimization code written in matlab.
  Structural and multidisciplinary optimization 21(2):120--127

\bibitem[{Sigmund(2007)}]{sigmund2007morphology}
Sigmund O (2007) Morphology-based black and white filters for topology
  optimization. Structural and Multidisciplinary Optimization 33(4):401--424

\bibitem[{Sigmund and Maute(2012)}]{sigmund2012sensitivity}
Sigmund O, Maute K (2012) Sensitivity filtering from a continuum mechanics
  perspective. Structural and Multidisciplinary Optimization 46(4):471--475

\bibitem[{Sigmund and Maute(2013)}]{sigmund_topoopt_approaches}
Sigmund O, Maute K (2013) Topology optimization approaches. Structural and
  Multidisciplinary Optimization 48(6):1031–--1055

\bibitem[{Str{\"o}mberg(2013)}]{stromberg2013influence}
Str{\"o}mberg N (2013) The influence of sliding friction on optimal topologies.
  In: Recent advances in contact mechanics. Springer, p 327--336

\bibitem[{Svanberg(1987)}]{mma_method}
Svanberg K (1987) The method of moving asymptotes—a new method for structural
  optimization. International journal for numerical methods in engineering
  24(2):359--373

\bibitem[{Wei{\ss}enfels and Wriggers(2015)}]{weissenfels2015contact}
Wei{\ss}enfels C, Wriggers P (2015) A contact layer element for large
  deformations. Computational Mechanics 55(5):873--885

\bibitem[{Wriggers et~al(2013)Wriggers, Schröder, and
  Schwarz}]{Wriggers_third_medium}
Wriggers P, Schröder J, Schwarz AA (2013) A finite element method for contact
  using a third medium. Computational Mechanics volume 52:837--847

\bibitem[{Zhang et~al(2021)Zhang, Li, Wang, Jia, and Luo}]{zhang2021narrow}
Zhang X, Li Y, Wang Y, et~al (2021) Narrow-band filter design of phononic
  crystals with periodic point defects via topology optimization. International
  Journal of Mechanical Sciences 212:106,829

\bibitem[{Zhou et~al(2016)Zhou, Zhang, Zhu, and Xu}]{zhou2016feature}
Zhou Y, Zhang W, Zhu J, et~al (2016) Feature-driven topology optimization
  method with signed distance function. Computer Methods in Applied Mechanics
  and Engineering 310:1--32

\end{thebibliography}

\end{document}